\documentclass[a4paper,10pt,final,leqno]{amsart}
\usepackage{amsthm,amsmath,amssymb,amscd,amsxtra,amsgen,enumitem,verbatim}
\usepackage{hyperref,ae}
\usepackage[?]{amsrefs}

\usepackage[all]{xy}
\newtheorem{thmlabel}{Theorem}

\ifx\undefined\xspace \usepackage{xspace} \fi

\newcommand{\nc}{\newcommand}
\renewcommand{\frak}{\mathfrak}
\providecommand{\cal}{\mathcal}
\renewcommand{\bold}{\mathbf}

\numberwithin{equation}{section}

\newcommand{\pfname}{Proof.}

\newenvironment{pf}{\vskip-\lastskip\vskip\medskipamount{\it\pfname}}%
                      {$\square$\vskip\medskipamount\par}

\newenvironment{pfof}[1]{\vskip-\lastskip\vskip\medskipamount{\it
    Proof of #1.}}%
                      {$\square$\vskip\medskipamount\par}

\newtheorem{thm}{Theorem}[subsection]

\newtheorem{thmnn}{Theorem}

\newtheorem{cor}[thm]{Corollary}

\newtheorem{prop}[thm]{Proposition}

\newtheorem{lemma}[thm]{Lemma} 

\theoremstyle{definition}
\newtheorem{defn}[thm]{Definition}
\newtheorem{definition}[thm]{Definition}

\theoremstyle{definition}

\newtheorem{remark}[thm]{Remark}
\newtheorem{remarks}[thm]{Remarks}

\newtheorem{example}[thm]{Example} 
\newtheorem{examples}[thm]{Examples} 

\nc{\Theorem}[1]{Theorem~{#1}}
\nc{\Th}[1]{({\sl Th.}~#1)}
\nc{\Thd}[2]{({\sl Th.}~{#1} {#2})}
\nc{\Theorems}[2]{Theorems~{#1} and ~{#2}}
\nc{\Thms}[2]{({\it Thms. ~{#1} and ~{#2}})}
\nc{\Lemmas}[2]{Lemma~{#1} and ~{#2}}

\nc{\manga}[6]{({\it Thms. ~ #1, ~ #2, ~ #3,\\ ~ #4, ~ #5, ~ #6})}

\nc{\Prop}[1]{({\sl Prop.}~{#1})}
\nc{\Proposition}[1]{Proposition~{#1}}
\nc{\Propositions}[2]{Propositions~{#1} and ~{#2}}
\nc{\Props}[2]{({\sl Props.}~{#1} and ~{#1})}
\nc{\Cor}[1]{({\sl Cor.}~{#1})}
\nc{\Corollary}[1]{Corollary~{#1}}
\nc{\Corollaries}[2]{Corollaries~{#1} and ~{#2}}
\nc{\Definition}[1]{Definition~{#1}}
\nc{\Defn}[1]{({\sl Def.}~{#1})}
\nc{\Lemma}[1]{Lemma~{#1}} 
\nc{\Lem}[1]{({\sl Lem.} ~{#1})} 
\nc{\Eq}[1]{equation~({#1})}
\nc{\Equation}[1]{Equation~({#1})}
\nc{\Section}[1]{Section~{#1}}
\nc{\Sections}[1]{Sections~{#1}}
\nc{\Sec}[1]{({\sl Sec.} ~{#1})} 
\nc{\Chapter}[1]{Chapter~{#1}}
\nc{\Chapt}[1]{({\sl Ch.}~{#1})}

\nc{\Ex}[1]{{\sl Ex.}~{#1}}
\nc{\Exa}[1]{{\sl Example}~{#1}}
\nc{\Example}[1]{{\sl Example}~{#1}}
\nc{\Examples}[1]{{\sl Examples}~{#1}}
\nc{\Exercise}[1]{{\sl Exercise}~{#1}}

\nc{\Rem}[1]{({\sl Rem.}~{#1})}
\nc{\Remark}[1]{{\sl Remark}~{#1}}
\nc{\Remarks}[1]{{\sl Remarks}~{#1}}
\nc{\Note}[1]{{\sl Note}~{#1}}

\nc{\Conjecture}[1]{Conjecture~{#1}}
\nc{\Claim}[1]{Claim~{#1}}

\nc \Proof{{  \it Proof. }}

\nc{\xmu}{\mu}
\nc{\w}{\omega}

\nc \Ab{{\ensuremath{\bold A}}}
\nc \ab{{\ensuremath{\bold a}}}
\nc \bb{{\ensuremath{\bold b}}}
\nc \cb{{\ensuremath{\bold c}}}
\nc \Bb{{\ensuremath{\bold B}}}
\nc \Gb{{\ensuremath{\bold G}}}
\nc \Qb{{\ensuremath{\bold Q}}}
\nc \Rb{{\ensuremath{\bold R}}} \nc \Cb{{\ensuremath{\bold C}}} 
\nc \Eb{{\ensuremath{\bold E}}}
\nc \eb{{\ensuremath{\bold e}}}
\nc \Db{{\ensuremath{\bold D}}}
\nc \Fb{{\ensuremath{\bold F}}}
\nc \ib{{\ensuremath{\bold i}}}
\nc \jb{{\ensuremath{\bold j}}}
\nc \kb{{\ensuremath{\bold k}}}
\nc \nb{{\ensuremath{\bold n}}}
\nc \rb{{\ensuremath{\bold r}}}
\nc \Pb{{\ensuremath{\bold P}}}
\nc \pb{{\ensuremath{\bold p}}}
\nc \SPb{{\ensuremath{\bold {SP}}}}
\nc \Zb{{\ensuremath{\bold Z}}} 
\nc \zb{{\ensuremath{\bold z}}} 
\nc \gb{{\ensuremath{\bold g}}} 
\nc \fb{{\ensuremath{\bold f}}} 
\nc \ub{{\ensuremath{\bold u}}} 
\nc \vb{{\ensuremath{\bold v}}} 
\nc \yb{{\ensuremath{\bold y}}} 
\nc \xb{{\ensuremath{\bold x}}} 
\nc \xib{{\ensuremath{\bold \xi}}} 
\nc \Nb{{\ensuremath{\bold N}}} 
\nc \Hb{{\ensuremath{\bold H}}} 
\nc \wb{{\ensuremath{\bold w}}} 
\nc \Wb{{\ensuremath{\bold W}}} 
\nc \syz{{\mathbf {syz}}}
\nc \bnoll{{\ensuremath{\bold 0}}} 

\nc \mf{\frak m} \nc \mh{\hat{\m}} 
\nc \nf{\frak n}
\nc \Of{\frak O}
\nc \rf{\frak r}
\nc \mufr{{\mathbf \mu}}
\nc \hf{\frak h} 
\nc \qf{\frak q} 
\nc \bfr{\frak b} 
\nc \kfr{\frak k} 
\nc \pfr{\frak p} 
\nc \af{\frak a }
\nc \cf{\frak c }
\nc \sfr{\frak s} 
\nc \ufr{\frak u} 
\nc \g{\frak g} 
\nc \gA{\g_{\Ao}} 
\nc \lfr{\frak l}
\nc \afr{\frak a}
\nc \gfh{\hat {\frak g}}
\nc \gl{\frak { gl }}
\nc \Sl{\frak {sl}}
\nc \SU{\frak {SU}}
\nc{\Homf}{\frak{Hom}}

\newcommand{\on}{\operatorname}
\nc\hankel{\on {Hankel}}
\nc\row{\on {row\ }}
\nc\nullity{\on {nullity }}
\nc\col{\on {col\ }}
\nc\rowm{\on {Row \ }}
\nc\loc{\on {lc \ }}
\nc\nullo{\on {null\ }}
\nc\Nul{\on {Nul\ }}
\nc \Ann {\on {Ann }}
\nc \Ass {\on {Ass \ }}
\nc \Coker {\on {Coker}}
\nc \Co{\on C}
\nc \Homo{\on {Hom}}
\nc \Ker {\on {Ker}}
\nc \omod{\on {mod}}
\nc \No {\on N}
\nc \NN {\on {NN}}
\nc \NGo {\on {NG}}
\nc \Oo {\on O}
\nc \ch {\on {ch}}
\nc \rko {\on {rk}}
\nc \Sing {\on {Sing\ }}
\nc \Reg {\on {Reg}}
\nc \CoI {\on {CI}}
\nc \CoM {\on {CM}}
\nc \Gor {\on {Gor}}
\nc \Type {\on {Type}}
\nc \can {\on {can}}
\nc \Top {\on {T}}
\nc \Tr {\on {Tr}}
\nc \rel {\on {rel}}
\nc \tr {\on {tr}}
\nc \sgn {\on {sgn }}
\nc \trdeg {\on {tr.deg}}
\nc \codim {\on {codim }}
\nc \coht {\on {coht}}
\nc \divo {\on {div \ }}
\nc \coh {\on {coh}}
\nc \Clo {\on {Cl}}
\nc \embdim{\on {embdim}}
\nc \embcodim{\on {embcodim \ }}
\nc \qcoh {\on {qcoh}}
\nc \grad {\on {grad}\ }
\nc \grade {\on {grade}}
\nc \hto {\on {ht}}
\nc \depth {\on {depth}}
\nc \prof {\on {prof}}
\nc \reso{\on {res}}
\nc \ind{\on {ind}}
\nc \prodo{\on {prod}}
\nc \coind{\on {coind}}
\nc \Con{\on {Con}}
\nc \Crit{\on {Crit}}
\nc \Der{\on {Der}}
\nc \Char{\on {Char}}
\nc \Ch{\on {Ch}}

\nc \Ext{\on {Ext}}
\nc \Eo{\on {E}}
\nc \End{\on {End}}
\nc \ad{\on {ad}}
\nc \Ad{\on {Ad}}
\nc \gr{\on {gr}}
\nc \Fo{\on {F}}
\nc \Gr{\on {Gr}}
\nc \Go{\on {G}}
\nc \GFo{\on {GF}}
\nc \Glo{\on {Gl}}
\nc \Ho{\on {H}}
\nc \CMo{\on {\CM}}
\nc \SCM{\on {SCM}}
\nc \hol{\on {hol}}
\nc{\sgd}{\on{sgd}}
\nc \supp{\on {supp}}
\nc \ssupp{\on {s-supp}}
\nc \singsupp{\on {singsupp}}
\nc \msupp{\on {msupp}}
\nc \spec{\on {spec}}
\nc \spano{\on {span }}
\nc \Span{\on {Span }}
\nc \Max{\on {Max}}
\nc \Min{\on {Min}}
\nc \Mod{\on {Mod}}
\nc \Rad {\on {Rad}}
\nc \rad {\on {rad}}
\nc \rank {\on {rank}}
\nc \range {\on {range}}
\nc \Slo{\on {SL}}
\nc \soc {\on {soc}}
\nc \Irr {\on {Irr}}
\nc \Imo {\on {Im}}
\nc \SSo{\on {SS}}
\nc \lub{\on {lub}}
\nc \gldim{\on {gl.d.}}
\nc \pdo{\on {p.d.}} 
\nc \ido{\on {i.d.}} 
\nc \dSSo{\dot {\SSo}}
\nc \So{\on S}
\nc \Io{\on I}
\nc \Jo{\on J}
\nc \jo{\on j}
\nc \Ko{\on K}
\nc \PBW{\Ac_{PBW}}
\nc \Ro{\on R}
\nc \To{\on T}
\nc \Ao{\on A}

\nc \Do{{\on D}}
\nc \Bo{\on B}
\nc \Po{\on P}
\nc \Qo{\on Q}
\nc \Zo{\on Z}
\nc \U{\on U}
\nc \wt{\on {wt}}
\nc \Uh{\hat {\U}}
\nc \T{\on T}
\nc \Lo{\on L}
\nc{\dop}{\on d}
\nc{\eo}{\on e}
\nc{\ado}{\on{ad}}
\nc{\Tot}{\on{Tot}}
\nc{\Aut}{\on{Aut}}
\nc{\sinc}{\on {sinc}}

\nc{\overrightleftarrows}[2]{\overset{#1}{\underset{#2}{\rightleftarrows}}}

\nc{\CCF}{\cal{CF}}
\nc{\CDF}{\cal{DF}}
\nc{\CHC}{\check{\cal C}}

\nc{\Cone}{\on{Cone}}
\nc{\dec}{\on{dec}}
\nc{\Diff}{\on{Diff}}
\nc{\dirlim}{\underset{\to}{\on{lim}}}
\nc{\dpar}{\partial}
\nc{\GL}{\on{GL}}
\nc{\CGr}{\cal{G}r}
\nc{\pr}{\on{pr}}
\nc{\semid}{|\!\!\!\times}
\nc{\Hom}{\on{Hom}}
\nc \RHom{\on {RHom}}

\nc \Proj{\mathrm {Proj\ }}
\nc \proj{\mathrm {proj}}
\nc{\Id}{\on{Id}}
\nc{\id}{\on{id}}
\nc{\Ima}{\on{Im}}
\nc{\invtimes}{\underset{\gets}{\otimes}}
\nc{\invlim}{\underset{\gets}{\on{lim}}}
\nc{\Lie}{\on{Lie}}
\nc{\re}{\on{Re }}
\nc{\Pic}{\on{Pic }}
\nc{\LPic}{\on{LPic }}
\nc{\Sch}{\on{Sch}}
\nc{\Sh}{\on{Sh}}
\nc{\Set}{\on{Set}}
\nc{\spo}{\on{sp\  }}
\nc{\Spec}{\on{Spec}}
\nc{\mSpec}{\on{mSpec}}
\nc{\Specb}{\bold {Spec}}
\nc{\Projb}{\bold {Proj}}
\nc{\Specan}{\on{Specan}}
\nc{\Spo}{\on{Sp}}
\nc{\Spf}{\on{Spf}}
\nc{\sym}{\on{sym}}
\nc{\symm}{\on{symm}}
\nc{\rop}{\on{r}}
\nc{\Td}{\on{Td}}
\nc{\Tor}{\on{Tor}}

\nc{\Artin}{\cal{A}rtin}
\nc{\Dgcoalg}{\cal{D}gcoalg}
\nc{\Dglie}{\cal{D}glie}
\nc{\Ens}{\cal{E}ns}
\nc{\Fsch}{\cal{F}sch}
\nc{\Groupoids}{\cal{G}roupoids}
\nc{\Holie}{\cal{H}olie}
\nc{\Mor}{\cal{M}or}

\nc{\CF}{\ensuremath{\cal{F}}}
\nc \Kc{\ensuremath{\cal K}}
\nc \Lc{{\ensuremath{\cal L}}}
\nc \lcc{{\mathcal l}} 
\nc \CC{{\ensuremath{\cal C}}} 
\nc \Cc{{\ensuremath {\cal C}}}
\nc \Pc{{\ensuremath{\cal P}}}
\nc \Dc{\ensuremath{\mathcal D}}
\nc \Ac{{\ensuremath{\cal A}}} 
\nc \Bc{{\ensuremath{\cal B}}}
\nc \Ec{{\ensuremath{\cal E}}}
\nc \Fc{{\ensuremath{\cal F}}}
\nc \Mcc{{\ensuremath{\cal M}}} 
\nc \hM{\hat{\Mcc}} 
\nc \bM{\bar {\Mcc}} 
\nc\hbM{\hat{\bar \Mcc}}  
\nc \Nc{{\ensuremath{\cal N}}}
\nc \Hc{{\ensuremath{\cal H}}} 
\nc \Ic{{\ensuremath{\cal I}}} 
\nc \Oc{\ensuremath{{\cal O}}}
\nc \Och{\hat{\cal O}} 
\nc \Sc{{\ensuremath{{\cal S}}}}
\nc \Tc{\ensuremath{{\cal T}}} 
\nc \Vc{{\ensuremath{{\cal V}}}} 
\nc{\CA}{{\ensuremath{{\cal A}}}}
\nc{\CB}{{\ensuremath{{\cal B}}}}
\nc{\C}{{\ensuremath{{\cal F}}}}
\nc{\Gc}{{\ensuremath{{\cal G}}}}
\nc{\CH}{\ensuremath{\mathcal H}}
\nc{\CI}{{\ensuremath{{\cal I}}}}
\nc{\CM}{{\ensuremath{{\cal M}}}}
\nc{\CN}{{\ensuremath{{\cal N}}}}
\nc{\CO}{{\ensuremath{{\cal O}}}}
\nc{\Rc}{{\ensuremath{{\cal R}}}}
\nc{\CT}{{\ensuremath{\mathcal T}}}
\nc{\CU}{\ensuremath{{\cal U}}}
\nc{\CV}{\ensuremath{{\cal V}}}
\nc{\CZ}{\ensuremath{{\cal Z}}}
\nc{\Homc}{\ensuremath{{\cal {Hom}}}}

\nc{\fa}{\frak{a}}
\nc{\fA}{\frak{A}}
\nc{\fg}{\frak{g}}
\nc{\fh}{\frak{h}}
\nc{\fI}{\frak{I}}
\nc{\fK}{\frak{K}}
\nc{\fm}{\frak{m}}
\nc{\fP}{\frak{P}}
\nc{\fS}{\frak{S}}
\nc{\ft}{\frak{t}}
\nc{\fX}{\frak{X}}
\nc{\fY}{\frak{Y}}

\nc{\bF}{\bar{F}}
\nc{\bCP}{\bar{\cal{P}}}
\nc{\bm}{\mbox{\bf{m}}}
\nc{\bT}{\mbox{\bf{T}}}
\nc{\hB}{\hat{B}}
\nc{\hC}{\hat{C}}
\nc{\hP}{\hat{P}}
\nc{\htest}{\hat P}

\nc{\nen}{\newenvironment}
\nc{\ol}{\overline}
\nc{\ul}{\underline}
\nc{\ra}{\to}
\nc{\lla}{\longleftarrow}
\nc{\lra}{\longrightarrow}
\nc{\Lra}{\Longrightarrow}
\nc{\Lla}{\Longleftarrow}
\nc{\Llra}{\Longleftrightarrow}
\nc{\hra}{\hookrightarrow}
\nc{\iso}{\overset{\sim}{\lra}}

\nc{\dsize}{\displaystyle}
\nc{\sst}{\scriptstyle}
\nc{\tsize}{\textstyle}
\nen{exa}[1]{\label{#1}{\bf Example.\ } }{}

\nen{rem}[1]{\label{#1}{\em Remark.\ } }{}
\nen{exer}[1]{\label{#1}{\em Exercise.\ } }{}

\begin{document}


\title[Liftable derivations]{Liftable derivations for generically separably algebraic
  morphisms of schemes} \author[\ ]{Rolf K{\"a}llstr{\"o}m}

\address{Department of Mathematics, University of G{\"a}vle, 801 76
  G{\"a}vle, Sweden} \email{rkm@hig.se}  \date \today
 \begin{abstract}
   We consider dominant, generically algebraic (e.g. generically
   finite), and tamely ramified (if the characteristic is positive)
   morphisms $\pi : X/S \to Y/S$ of $S$-schemes, where $Y,S$ are
   N{\oe}therian and integral and $X$ is a Krull scheme (e.g. normal
   N{\oe}therian), and study the sheaf of tangent vector fields on $Y$
   that lift to tangent vector fields on $X$. We give an easily
   computable description of these vector fields using valuations
   along the critical locus. We apply this to answer the question when
   the liftable derivations can be defined by a tangency condition
   along the discriminant.  In particular, if $\pi$ is a blow-up of a
   coherent ideal $I$ we show that tangent vector fields that preserve
   the Ratliff-Rush ideal (equals $[I^{n+1}:I^n]$ for high $n$)
   associated to $I$ are liftable, and that all liftable tangent
   vector fields preserve the integral closure of $I$. We also
   generalise to positive characteristic Seidenberg's theorem that all
   tangent vector fields can be lifted to the normalisation, assuming
   tame ramification.
 \end{abstract}

\maketitle
\noindent
\section*{Introduction} 
Let $\pi: X/S\to Y/S$ be a dominant morphism of $S$-schemes, where $Y$
and $S$ are N{\oe}therian and $X$ is a Krull scheme (e.g.  N\oe therian
and normal), so in particular $X$ and $Y$ are integral.  Consider the
diagram
\begin{displaymath}
 \xymatrix{T_{X/S}\ar[r]^{d\pi\ \ \ }& T_{X/S \to Y/S}\ar[r]& \Cc_{X/Y}\ar[r]& 0\\
& \pi^{-1}(T_{X/S})\ar[u]_{\psi}   & &}
\end{displaymath}
where $T_{X/S}=Hom_{\Oc_X}(\Omega_{X/S}, \Oc_X)$ is the sheaf of
$S$-relative tangent vector fields, $d\pi$ is the tangent morphism and
$\psi $ is the canonical morphism to the sheaf $T_{X/S \to Y/S}$ of
derivations from $\pi^{-1}(\Oc_Y)$ to $\Oc_X$; if $Y/S$ is locally of
finite type and either $\pi$ is flat or $Y/S$ is smooth, then $ T_{X/S
  \to Y/S}=\pi^*(T_{Y/S})$. The critical locus $C_\pi$ is the support of
the critical module $\Cc_{X/Y}$, and the discriminant set $D_\pi$ is
the closure of $\pi(C_\pi)$.  Say that a section $\partial$ of $T_{Y/S}$ is
liftable if there exists a section $\bar \partial$ of $T_{X/S}$ such that
$d\pi(\bar \partial ) = \psi (\partial)$ (suppressing domains of definition of $\partial$
and $\bar \partial$).  We denote the sheaf of liftable tangent vector fields
by $T^\pi_{Y/S}$.

Let $T_{Y/S}(\log^\pi I_{D_\pi}) = T_{Y/S}(\log^\pi D_\pi)$ be the sheaf
of {\it $\pi$-logarithmic vector fields}.  This sheaf is easily
computed using the discrete valuations $\nu_1, \dots , \nu_r$ of the
function field of $X$ defined by order of vanishing along the
components of the critical divisor locus $C_\pi = C_1 + C_2 + \cdots +
C_r$. A section $\partial $ in $T_{Y/S}$ is $\pi$-logarithmic, or logarithmic
along the morphism $\pi$, if it does not lower the value for the
discrete valuations at the generic points of $C_\pi$ of local
generators of the ideal $I_{D_\pi}$, at points below the generic points,
i.e.  $\nu_i(\partial (f_j))\geq \nu_i(f_j)$ when $f_j \in I_i$, where $I_i$ is
any defining ideal of the discriminant $D_{\pi}$ at the generic point
of $\pi(C_i)$. The sheaf of {\it logarithmic vector fields} along a
subscheme defined by an ideal $I$, denoted $T_{Y/S}(I)$, is the sheaf
of tangent vector fields $\partial$ that preserve the ideal $I$, $\partial(I)\subset
I$. We describe the liftable tangent vector fields when $\pi$ is
generically algebraic (the induced extension of function fields is
algebraic), and tamely ramified \Defn{\ref{tame-ramification}}.

\begin{thmlabel} Let $\pi: X/S \to Y/S$ be a morphism of schemes over a
  scheme $S$.  Assume that $Y$ and $S$ are N{\oe}therian and $X$ is a
  Krull scheme.  Assume that $\pi$ is dominant, generically algebraic
  and tamely ramified. If $Y$ is of unequal characteristic, assume
  that $Y$ is absolutely tamely ramified. Let $D_\pi$ be the discriminant
  set of $\pi$. Then:
    \begin{enumerate}[label=($L_\theenumi$),ref=($L_\theenumi$)]\label{main-intro}
    \item Let $J_{D_\pi}$ be {\it any}\/ defining ideal of $D_\pi$.
      Then
      \begin{displaymath} T_{Y/S}^\pi = T_{Y/S}(\log^\pi J_{D_\pi}).
      \end{displaymath}
      \label{L1i}\\
    \item If the residue field extension $k_c/k_{\pi(c)}$ is algebraic
      for each point $c\in C_\pi$ of height $1$, then
      \begin{displaymath} T_{Y/S}^\pi = T_{Y/S}(I_{D_\pi}),\label{L2i}
      \end{displaymath}
      where $I_{D_\pi}$ is the (reduced) ideal of the closure of the
      discriminant set $D_\pi$.
    \item Let $I_{C_\pi}$ be the ideal of the critical locus of $\pi$.
      Then $\pi^{-1}(T_{Y/S}^\pi) \subset T_{X/S}(I_{C_\pi})$, i.e. the
      liftable vector fields are tangent to the critical
      locus.\label{L3i}
    \end{enumerate}
  \end{thmlabel}
  We give concrete examples at the end of the paper.

  Here \ref{L1i} is the main new content of \Theorem{\ref{main-intro}}
  while \ref{L2i} is a generalisation of earlier work, but then
  usually formulated for affine schemes.  The first result of this
  type is due to Zariski who proved the following
  \cite{zariski:equisinglarI}*{{\S} 5, Th.  2}.  Let $R$ be an
  integrally closed local Nagata ring, with quotient field $K(R)$.
  Let $R'$ be the integral closure in a finite separable extension $L$
  of $K(R)$ (thus $R'$ is finite over $R$).  Let $D\subset \Spec R$ be the
  set of ramified primes of height $1$ in $R$, and assume that each
  prime in $D$ is tamely ramified.  Let $\partial$ be a derivation of $R$
  such that for each point $x\in D$ of height $1$ there exists an
  element $r\in \mf_x\setminus \mf_x^ 2$ such that $\partial(r)=0$.  It then follows
  that $\partial (R')\subset R'$.  According to \Theorem{\ref{mainthm}}, with $X=
  \Spec R'$ and $Y=\Spec R$, knowing the existence of the element $r$
  is not necessary, as is the assumption that $R$ be a Nagata ring;
  the sufficient and necessary condition is~\ref{L2i}.  Scheja and
  Storch \cite{scheja:fortzetsungderivationen} proved~\ref{L2i} when
  $R$ contains the rational numbers, assuming the $R$-module $R\partial (R)$
  be of finite type.
We see that this assumption is not
needed and we also get the assertion in positive characteristic when
$R/A$ is tamely ramified.

A dominant morphism $\pi$ is {\it weakly submersive}\/ if $T^\pi_{Y/S}=
T_{Y/S}$, i.e.\ all tangent vector fields lift to vector fields on
$X$. If $Y/S$ is non-smooth and $\pi$ is non-flat this is {\it not} the
same as surjective tangent morphism $d\pi$, in which case we say $\pi$
is submersive.  Seidenberg
\cite{seidenberg:derivationsintegralclosure} proved that the
normalisation morphism of a N{\oe}therian integral scheme $X/S$ is
weakly submersive when $S$ is defined over the rational numbers.  We
generalise this to positive characteristic \Th{\ref{seidenberg}},
requiring that the normalisation morphism be tamely ramified. The
proof is by reducing to the case when the normalisation is a discrete
valuation ring, which is quite different from Seidenberg's proof (see
\Remark{\ref{remark-seidenberg}}). In a forthcoming paper we will
prove that the constructive resolutions of singularities of a variety
$X/k$ ($\Char k =0$) presented in
\cites{villamayor:res,bierstone-milman:resolution,encinas-villamayor:good}
are weakly submersive. This also gives the result that the multiplier
ideal $J(\alpha)$ (discussed in \cite{lazarsfeld:vol2}) of a coherent
ideal $I$ is preserved by derivations that preserve $I$, i.e.
$T_{X/k}(I)\subset T_{X/k}(J(\alpha))$.

\Theorem{\ref{mainthm}} describes the liftable derivations at generic
points of the discriminant $D_\pi$ in terms of a set of discrete
valuations associated to the critical locus.  A natural question is
whether the liftable tangent fields coincide with the vector fields
that are tangent to some subscheme whose underlying space is $D_\pi$.
More precisely, does there exist a coherent defining ideal $I$ of
$D_\pi$ satisfying
\begin{equation}\tag{$*$}
  T^\pi_{Y/S}= T_{Y/S}(I)?
\end{equation}
The best one can hope for, second to weakly submersive, is that all
vector fields that are tangent to the discriminant can be lifted, i.e.
the above equality holds when $I$ is the (radical) ideal of $D_\pi$. We
call such morphisms {\it differentially ramified}, and know already
from \Theorem{\ref{mainthm}} that residually algebraic tamely ramified
morphisms are differentially ramified.  A stronger notion is that $\pi$
be {\it uniformly ramified}, meaning that at each generic point $\xi $
of the critical locus $C_\pi$, the stalk $I_{D_\pi,\pi(\xi)}$ has a basis
with constant value in the valuation ring $\Oc_{X,\xi}$.  We prove in
\Theorem{\ref{blowup1}} that uniformly ramified morphisms are
differentially ramified, and if for each generic point $\xi\in C_\pi$ there
exists a basis $\{x_1, \dots , x_r\}$ of $I_{D,\pi(\xi)}$ and
derivations $\{\partial_1, \dots , \partial_r\} \subset T_{Y/S, \pi(\xi)}$ such that the
matrix $(\partial_j(x_i))$ is invertible, then the converse also holds.

It is particularly interesting to characterise the liftable tangent
fields for birational morphisms, forming an important class of
generically separably algebraic morphisms. In the light of \ref{L2i}
and \ref{L3i} in \Theorem{\ref{mainthm}} this also characterises the
liftable vector fields for ``alterations'', i.e.  compositions of
finite and birational morphisms. We have already discussed certain
birational morphisms, but now consider any (projective) birational
morphism of integral N{\oe}therian schemes, which we know always is
the blow-up of some fractional ideal. Hence let $\pi: Bl_I(X)\to X$ be
the blow-up of a given fractional ideal $I$ on $X$.  Here $Bl_I(X)$
need not be Krull (normal) so \Theorem{\ref{mainthm}} is not directly
applicable, but we can apply it to attain a lower and upper inclusion
for $T^\pi_{X/S}$ by sheaves of vector fields that are tangent to
certain subschemes with the same underlying space as the discriminant.
One first easily gets that the sub-Lie algebroid of vector fields that
preserve $I$ are liftable, $T_{X/S}(I)\subset T^\pi_{X/S}$, but in
general this inclusion is strict.  There is also a latitude in the
choice of $I$, for different ideals may give the same blow-up, and we
ask if there is a choice $\tilde I$ such that $Bl_I(X)\cong Bl_{\tilde
  I}(X)$ (isomorphism over $X$) and $T^\pi_{X/S} = T_{X/S}(\tilde I)$?
Let $\hat I =\cup_{n\geq 1} [I^{n+1}: I^n]$ be the Ratliff-Rush ideal
associated to $I$ \cite{ratliff-rush} and $\bar I $ its integral
closure; then $I \subset \hat I \subset \bar I$. In
\Theorem{\ref{blow-up-thm}} we prove
\begin{displaymath}
  T_{X/S}(\hat I) \subset T^\pi_{X/S} \subset T_{X/S}(\bar I),
\end{displaymath}
where the right-hand inclusion holds under the additional assumption
that the normalisation morphism of $Bl_I(X)$ be tamely ramified (make
this assumption in this paragraph). We therefore succeed in a positive
answer to our question $(*)$ above when the Ratliff-Rush ideal $\hat I$ associated
to $I$ is integrally closed, getting $T^\pi_{X/S} = T_{X/S}(\hat I)$,
since $Bl_{\hat I}(X)$ is isomorphic to $Bl_I(X)$.  An immediate
consequence is that blow-ups of radical ideals are differentially
ramified.  Combining with \Theorem{\ref{blowup1}} one also gets this:
Let $\pi : \tilde X/S \to X/S$ be a tamely ramified blow-up of a reduced
subscheme $V/S$ where $\tilde X$ is assumed to be Krull and $X/S$ is a
N{\oe}therian integral scheme which is smooth at the generic points of
$V$, then $\pi$ is uniformly ramified \Cor{\ref{cor:blow}}.

One should note that several results in this paper have
straightforward holomorphic counterparts.  For instance, for a finite
holomorphic map $\pi : \Cb^m \to \Cb^m$ we have $T^\pi_{\Cb^m}=
T_{\Cb^m}(I_{D_\pi})$, which was proven in a different way in
\cite{arnold:wave}.

We also mention how liftable tangent vector fields can be used.  One
may interprete $T^\pi_{Y/S}$ is as the sheaf of infinitesimal
symmetries of $\pi$, so that the fibres of liftable tangent vector
fields correspond to directions in the base where the fibres of $\pi$
do not deform. Thus a good understanding of $T^\pi_{Y/S}$ is useful for
the study of deformations in a (flat) family of schemes.  We tie this
up with Zariski's notion of analytic equisingularity stratification of
a hypersurface $X$ in $\Cb^n$.  The dominant stratum is where $X$ is
smooth, so that smooth points do not belong to the critical locus of
the restriction $\pi = p|_X : X \to \Cb^{n-1}$ for {\it generic}
projections $p: \Cb^n \to \Cb^{n-1}$.  Hence all tangent vector fields
near $\pi(x)$ lift to vector fields near a smooth point $x$; therefore,
by Zariski's lemma \cite[Corollary to Th.  30.1]{matsumura}, all
smooth points on $X$ can be regarded equisingular, since they have
isomorphic analytic neighbourhoods.  Of course, this must be so since
$\pi$ is {\'e}tale at $x$, so $\Oc_x$ is analytically isomorphic to the
ring or formal power series (or the analytic localisation of $\Oc_x$
is isomorphic to the ring of convergent power series in $n-1$
variables).  The next stratum consists of points $x$ where the
discriminant of $\pi_x$ is non-empty but smooth for generic $p$, so
that, at $\pi(x)\in D_\pi$, the stalk $T_{D_\pi,\pi(x)}$ of tangent vector
fields on $D_\pi$ near $\pi(x)$ has a basis of non-vanishing tangent
vector fields, which extend to non-vanishing vector fields near
$\pi(x)$ in $\Cb^{n-1}$ tangent to $D_\pi$, hence they are liftable.
Again by Zariski's lemma, points with smooth discriminants (and
generic $p$) can be regarded equisingular, since they have isomorphic
analytic neighbourhoods.  If the generic discriminants $D_\pi$ are not
smooth, by induction in the dimension $n$ the equisingularity
stratification of $D_\pi$ is defined, which can be pulled back to get a
stratification of $X$. For a more complete discussion of Zariski's
equisingularity stratification see
\cite{lipman:equisingular,villamayor:equisingular}.  Another situation
where a good description of $T^\pi_{Y/S}$ is useful, is for describing
direct images of D-modules (see \cite{borel:Dmod}) with respect to a
morphism $\pi$ of complex algebraic manifolds, say the zeroth direct
image of the structure sheaf $R^0 \pi_+(\Oc_X)=
\pi_*(\omega_{X/Y}\otimes_{\Dc_X}\pi^*(\Dc_Y))$ (in the complement of $D_\pi$ it
is a Gauss-Manin connection), where $\Dc_X$ and $\Dc_Y$ are the rings
of differential operators on $X$ and $Y$, and $\omega_{X/Y}$ is the
relative canoncial bundle. Then we have a surjective homomorphism of
$\Dc_Y$-modules
\begin{displaymath}
  \frac {\Dc_Y}{\Dc_Y T^\pi_{Y}}\to R^0 \pi_+(\Oc_X) \to 0,
\end{displaymath}
which in particular gives a bound on the singular support of
$R^0\pi_+(\Oc_X)$.

\bigskip {\it Standard notions}. Let $A$ and $R$ be local rings $(A,
\mf_A, k_A)$ and $(R, \mf_R, k_R)$, where $k_A, k_R$ are the residue
fields, and $\pi : A \to R$ a homomorphism of rings.  Then $\pi$ is local
if $\pi(\mf_A) \subset \mf_R$, $R$ dominates $A$ if $\pi^{-1}(\mf_R)= \mf_A$,
and $\pi$ is (separably) algebraic if the extension of fraction fields
$K(R)/K(A)$ is (separably) algebraic, residually
algebraic/finite/separable if $k_R/k_A$ is algebraic/finite/separable;
in particular, $R/A$ is residually finite if it is quasi-finite, while
the converse holds if $R\mf_A$ is an ideal of definition of $R$. We
refer to \cite{EGA4:1, matsumura} for the basic results about formally
smooth/unramified morphisms; the discrete topology is intended if no
other topology is mentioned. A morphism of schemes $\pi: X\to Y$ is
dominant if the morphism $\pi^{-1}(\Oc_Y)\to \Oc_X$ is injective; $\pi$
is generically (separably) algebraic if $\pi$ is dominant, $X$ and $Y$
are reduced, and for points $x\in X$ that map to a generic point $\xi \in
Y$ the corresponding residue field extension $k_x/k_{\xi}$ is
(separably) algebraic, and in particular formally \'etale. As common
practice, by $\Oc_X$-module we mean a sheaf of modules, and by writing
$m\in M$ for a sheaf $M$ we mean that $m$ is a section of $M$ over a
suitably defined open set. By a {\it Lie algebroid}\/ on $X/S$ is
intended an $\Oc_X$-module $\g_{X/S}$ which moreover is a Lie algebra,
provided with a homomorphism (as Lie algebras and $\Oc_X$-modules) to
the tangent sheaf $\alpha: \g_{X/S}\to T_{X/S}$, with the obvious
compatibility relations $[\delta,f\eta]= \alpha(\delta)(f)\eta + f[\delta, \eta]$, $\delta, \eta
\in \g_{X/S}$, $f\in \Oc_X$ (see e.g.  \cite{kallstrom:preserve}).  We
shall only have occasions to study Lie algebroids where $\alpha$ is
injective. For example, the sub-sheaf of liftable derivations
$T^\pi_{Y/S}$ is a sub-Lie algebroid of $T_{Y/S}$; other sub-Lie
algebroids arise from (fractional) ideals $I$ of $\Oc_Y$, as the
subsheaf $T_{Y/S}(I)$ of $T_{Y/S}$ of derivations $\partial$ that preserve
$I$, $\partial(I)\subset I$.

\section{Tangent morphisms, logarithmic derivations and
  ramification}\label{tangentmap}
The purpose of the material in \Sections{~\ref{quasi-finite-section}
  and~\ref{ram-section}}, which is not the minimal necessary to prove
our main results, is to provide a useful collection of basic
algebraicity/finiteness facts and some ramification theory, with
extensions of well-known results, and relate these classical notions
to our lifting problem.  For instance, we recall concisely the
relation between the property that a morphism $A\to R$ be finite and
quasi-finite, respectively, \Prop{\ref{zariski-main}} (a formulation
of Zariski's main theorem), and the relation between the property that
$K(R)/K(A)$ and $k_R/k_A$ be algebraic, respectively,
\Prop{\ref{q-finite}}.  We say a homomorphism of local N{\oe}therian
rings $ A \to R$ is ramified if the cotangent mapping $k_R\otimes \mf_A/
\mf_A^2 \to \mf_R/ \mf_R^2$ is not injective. This notion of
ramification is adequate for the study of derivations, and in the
important case when $R$ is a discrete valuation ring this is close to,
but stronger than, formally unramified; we describe what ramification
means for the liftable derivations \Prop{\ref{basic-prop}}. Formally
unramified/\'etale morphisms are described in some detail, complementing
the literature: \Proposition{\ref{etale-field}} is a characterisation
of separably {\it algebraic}\/ field extensions in terms of vanishing
differentials, and \Theorem{\ref{formunramified}} characterises
formally unramified morphisms of N{\oe}therian rings $R/A$ as morphisms
such that $k_R/k_A$ is separably algebraic and $R\mf_A = \mf_R$,
removing a finiteness assumption in Auslander and Buchsbaum's proof of
the assertion \cite{auslander-buchsbaum}.  \Lemma{\ref{p-rings}}
states that separably algebraic field extensions $k'/k$ extend to
$\mf$-\'etale extensions of complete $p$-rings (although not \'etale in
the discrete topology); this is needed in the proof of
\Theorem{\ref{mainthm}} to handle unequal characteristics.

\subsection{The tangent morphism}\label{first-sec} A morphism of
$S$-schemes $\pi : X/S \to Y/S$ has its canonical exact sequence of
differentials
\begin{equation}\label{eq:can-ex} \pi^*(\Omega_{Y/S}) \to \Omega_{X/S}
  \to \Omega_{X/Y} \to 0
\end{equation} which induces a homomorphism of
$\Oc_Y$-modules, the {\it tangent morphism}\/ of $\pi$,
\begin{displaymath} d\pi: T_{X/S}=Hom_{\Oc_X}(\Omega_{X/S},
\Oc_X) \to T_{X/S\to Y/S},
\end{displaymath} where the $\Oc_X$-module of `derivations
from $\Oc_Y$ to $\Oc_X$' is
\begin{displaymath} T_{X/S\to Y/S}
=Hom_{\Oc_X}(\pi^*(\Omega_{Y/S}),\Oc_X)=Hom_{\pi
^{-1}(\Oc_Y)}(\pi ^{-1}(\Omega_{Y/S}),\Oc_X).
\end{displaymath} There are canonical homomorphisms:
\begin{eqnarray}\label{can-tan0}
\psi_0&:&\pi^{-1}(T_{Y/S})\to T_{X/S\to Y/S},\\ \psi &:&
\pi^*(T_{Y/S}) \to T_{X/S\to Y/S}. \label{can-tan}
\end{eqnarray} If $\pi$ is dominant, then
$\psi_0$ is injective, but $\psi$ need be neither surjective
nor injective when $Y/S$ is non-smooth.

\begin{example}\label{ex:psi} (i) Let $X=\Spec k[t]$, $Y=\Spec
  k[t^2,t^3]$, and $S=\Spec k$.  Then $T_{X/k}= \Oc_X\partial_t$ and $T_{Y/k}= \Oc_Y
  t\nabla + \Oc_Y \nabla$ where $\nabla = t\partial_t$.  Let $\pi$ be
  the morphism of schemes induced by the inclusion of rings.  Then
  $\partial_t\in T_{X/S}$, and $t\otimes \nabla -1\otimes t\nabla \in
  \pi^*(T_Y)$ is non-zero, while $\psi (t\otimes \nabla -1\otimes
  t\nabla) = 0$. Also $d\pi(\partial_t)$ is the non-zero section of
  $T_{X/S\to Y/S}$ that is induced by the derivation $\partial_t:
  k[t^2,t^3] \to k[t]$; hence $d\pi(\partial_t)\notin \Imo (\psi)=
  \Oc_Y\nabla$. So $\psi$ is neither injective nor surjective. Here
  $d\pi$ is injective.

  (ii) (Immersion of singular locus) Assume that $X/k$ is a variety of
  characteristic $0$, with singular locus $Z$; let $I_{Z}$ be the
  ideal of $Z$. Let $\pi : Z \to X$ of the inclusion of the locus of
  non-smooth points and $T_{X/k}(\Oc_X,I_{Z})$ the subsheaf of
  $T_{X/k}$ consisting of derivations that send $\Oc_X$ to $I_{Z}$.
  We have an exact sequence
  \begin{displaymath} 0\to \pi^*(T_{X/k}(\Oc_X,I_{Z}))\to
\pi^*(T_{X/k})\xrightarrow{\psi } T_{Z/k \to X/k} \to
\Fc_{Z/X}\to 0
\end{displaymath} where the cokernel $\Fc_{Z/X}$ of
$\psi$ is a subsheaf of $Ext^1_{\Oc_X}(\Omega_{X/k},
I_{Z})$. The support of $ \Fc_{Z/X}$ consists of points
where $X$ is not equisingular, in the strong sense that
there exist tangent vector fields on $Z$ that cannot be extended  to tangent vector fields on $X$. Since $T_{X/k}=T_{X/k}(I_{Z})$ (see
e.g.  \cite{kallstrom:preserve}) there exists a natural
homomorphism $\can : \pi^*(T_{X/k})\to T_{Z/k}$, and a
commutative diagram:
  \begin{displaymath} \xymatrix{\pi^*(T_{X/k})
\ar[rd]^{\can}\ar[r]^{\psi}& T_{Z/k \to X/k}\\
&T_{Z/k}\ar[u]_{d\pi} }
  \end{displaymath} This gives a surjective homomorphism
$\Fc_{Z/X}\twoheadrightarrow \Cc_{Z/X}$. Therefore
$d\pi$ is surjective if $\psi$ is surjective.
\end{example}

\begin{prop}\label{psi} Assume that $Y/S$ is locally of
finite type.  The homomorphism $\psi$ is an isomorphism in
the following cases:
  \begin{enumerate}[label=(\theenumi)]
  \item $Y/S$ is smooth.
  \item $\pi$ is flat.
  \end{enumerate}
\end{prop}
\begin{pf} The map $\psi$ is an isomorphism if it induces an
isomorphism of stalks, and since $\Omega_{Y/S}$ is coherent
and $\Omega_{Y/S, y} \cong \Omega_{\Oc_{Y,y}/\Oc_{S,s}}$ (if
$U\subset Y$ is affine and $y\in U$, then the canonical
mapping $\Oc_{Y}(U)\to \Oc_{Y,y}$ is \'etale; then apply the
first fundamental exact sequence for differentials)
  \begin{displaymath} Hom_{\Oc_Y}(\Omega_{Y/S},\Oc_Y)_y =
Hom_{\Oc_{Y,y}}(\Omega_{\Oc_{Y,y}/\Oc_{S,s}}, \Oc_{Y,y}),
  \end{displaymath} where $s$ is the image of $y$ in $S$,
the condition is that the canonical morphism
  \begin{displaymath} \psi_x: \Oc_x \otimes _{\Oc_y}
\Hom_{\Oc_y}(\Omega_{\Oc_y/\Oc_s}, \Oc_y) \to
Hom_{\Oc_x}(\Oc_x \otimes_{\Oc_y} \Omega_{\Oc_y/\Oc_s},
\Oc_x)
  \end{displaymath} be an isomorphism, where $\Oc_y \to
\Oc_x$ is the homomorphism of local rings defined by $\pi$
and points $x\in X$, $y= \pi (x)\in Y$.  By (1)
$\Omega_{\Oc_y/\Oc_s}$ is free of finite rank over $\Oc_y$
\cite[Prop.  17.2.3]{EGA4:4}, implying the assertion.  Since
$\Omega_{\Oc_y/\Oc_s}$ is of finite presentation, (2)
implies the assertion by \cite[Th.  7.11]{matsumura}.
\end{pf}

When $\psi$ is an isomorphism one gets the ``ordinary''
tangent homomorphism
\begin{displaymath} \psi^{-1}\circ d\pi : T_{X/S} \to
\pi^*(T_{Y/S}).
\end{displaymath} If $\psi $ is not an isomorphism one can
take the fibre product to get a restricted tangent sheaf
\begin{multline*} T^r_{X/S}= \pi^*(T_{Y/S})\times_{T_{X/S\to Y/S}}T_{X/S} \\ =
  \{(\delta , \partial ) \in \pi^*(T_{X/S})\times _{T_{X/S\to Y/S}}T_{X/S} : \psi (\delta )= d\pi (\partial )
  \},
\end{multline*} so the projection on the first factor $T^r_{X/S} \to
\pi^*(T_{Y/S})$ can also play the role of `tangent morphism'.  The sheaf
$T^r_{X/S}$ is an $\Oc_Y$-Lie algebroid (namely the pull-back of the
Lie algebroid $T_{Y/S}$) containing the Lie sub-algebroid $T_{X/Y}=
\Ker (d\pi)$ of relative tangent vector fields and also the $\Oc_X$-{\it
  Lie algebra} $\bfr_X = \Ker (\psi)$. One gets a commutative diagram of
$\Oc_X$-modules:
\begin{equation}\label{diagram} \xymatrix{&&0&0&0\\
&0\ar[r]&N^4_X\ar[u]\ar[r]& \Cc_{X/Y}
\ar[u]\ar[r]&N^3_X\ar[u]\ar[r]&0\\
0\ar[r]&{\bfr}_X\ar[r]&\pi^*(T_{Y/S})\ar[r]^{\psi}\ar[u]&T_{X/S\to
Y/S}\ar[u]\ar[r]&N^2_X\ar[u]\ar[r]&0\\
0\ar[r]&{\bfr}_X\ar[r]\ar@{=}[u]&T^r_{X/S}\ar[r]\ar[u]&T_{X/S}\ar[u]_{d\pi}\ar[r]&N^1_X\ar[u]_h\ar[r]&0\\
&&T_{X/Y}\ar[u]\ar@{=}[r]&T_{X/Y}\ar[u]&0\ar[u]\\
&&0\ar[u]&0\ar[u]&&& }
\end{equation}

\subsection{Submersions and weak submersions}
\label{liftingderivations} Consider the vertical sequence above,
containing $d\pi$.  The stalk at $x\in X$ of the tangent morphism
\begin{displaymath} d \pi_x : T_{X/S,x} \to T_{X/S\to
Y/S,x}
\end{displaymath} induces a map of fibres $d \bar \pi_x:
k_x\otimes_{\Oc_x}T_{X/S,x} \to k_x\otimes_{\Oc_x}T_{X/S\to
  Y/S,x}$; if $T_{X/S\to Y/S}\cong \pi^*(T_{X/S})$
\Prop{\ref{psi}} one gets the ``ordinary'' map of tangent
spaces $k_x\otimes_{\Oc_x}T_{X/S,x} \to
k_x\otimes_{\Oc_y}T_{Y/S,y}$. Say that $\pi$ is {\it
  submersive}\/ at $x$ if $d\pi_x$ is surjective, i.e.
$\Cc_{X/Y,x}=0$.  If $\Cc_{X/Y,x}$ is of finite type, by
Nakayama's lemma $d \pi_x$ is surjective if the map $d
\bar \pi_x $ is surjective. The relation between the critical locus $C_\pi = \supp \Cc_{X/Y}$ and the ramification locus of a morphism is studied in \cite{kallstrom:purity}.


There are canonical homomorphisms of $\Oc_Y$-modules
$\pi_*(d \pi): \pi_*(T_{X/S})\to \\  \pi_*(T_{X/S\to Y/S})$
and $\pi_*(\psi_0):T_{Y/S}\to \pi_*(T_{X/S\to Y/S})$ so one
can form the fibre product
\begin{displaymath} \pi_*(T_{X/S})\times_{\pi_*(T_{X/S\to
Y/S})}T_{Y/S} = \{(\delta ,\partial) \in
\pi_*(T_{X/S})\times T_{Y/S} : \pi_*(d \pi )(\delta ) =
i(\partial )\}.
\end{displaymath}
\begin{defn} The sub-sheaf $T^\pi_{Y/S} \subset T_{Y/S}$ of
{\it liftable derivations}\/ on $Y$ is the image of the
projection on the second factor
  \begin{displaymath} T^\pi_{Y/S}= \Imo
(\pi_*(T_{X/S})\times_{\pi_*(T_{X/S\to Y/S})}T_{Y/S} \to
T_{Y/S}).
  \end{displaymath}
\end{defn} The {\it discriminant module} ${\Ec}_{X/Y}$ is the cokernel
of the inclusion morphism $T^\pi_{Y/S}\hookrightarrow T_{Y/S}$ and
$D_\pi^{disc} = \supp \Ec_{X/Y} $ the {\it weak discriminant set}\/ of
$\pi$.  Say that $\pi$ is {\it weakly submersive}\/ if $\Ec_{X/Y}=0$.
There are adjoint homomorphisms
\begin{eqnarray*} \phi_0 &:& \pi^{-1}(\Ec_{X/Y})\to
\Cc_{X/Y},\\ \phi&:& \Ec_{X/Y}\to \pi_*(\Cc_{X/Y}).
\end{eqnarray*}
\begin{prop}
    \begin{enumerate}[label=(\theenumi)]
  \item If $T_{X/Y}=0$ and $\psi_0$ (see~\ref{can-tan0}) is
injective, then $\phi_0$ is injective.
  \item If $R^1\pi_*(T_{X/Y})=0$ and the composed morphism
    \begin{displaymath} T_{Y/S}\to
\pi_*\pi^{-1}(T_{Y/S})\xrightarrow{\pi_*(\psi_0)}
\pi_*(T_{X/S\to Y/S})
  \end{displaymath} is injective, then $\phi$ is injective.
  \item If $\pi_*(T_{X/Y})=R^1\pi_*(T_{X/Y})=0$
(e.g. $T_{X/Y}=0$) and $\psi$ (see~\ref{can-tan}) is an
isomorphism \Prop{\ref{psi}}, hence $\pi_*(\psi)$ is an
isomorphism, then $T^\pi_{Y/S}= \pi_*(T_{X/S})$,
${\Ec}_{X/Y} = \pi_*(\Cc_{X/Y})$, and in particular $D_\pi^{disc}
=\pi( C_\pi)$.
  \end{enumerate} If either of the conditions in $(1)$ or
$(2)$ holds and $\pi$ is submersive, then $\pi$ is weakly
submersive.
\end{prop} Assuming the conditions in $(1)$, which is the situation we
shall mostly deal with, then in general when $\psi$ is not an
isomorphism we have a strict inclusion $D^{disc}_\pi \subset
\pi(C_\pi)$, so the notion of submersive morphism is stronger than
that of weakly submersive morphism; it is in fact straightforward to
see that these notions coincide (assuming (1)) if and only if $\psi$
is an isomorphism.  In \Example{\ref{ex:psi1}} below the morphism
$\pi$ is weakly submersive but not submersive; the fact that this
normalisation morphism is weakly submersive is not coincidental
\Th{\ref{seidenberg}}.

To describe the sheaf $T^\pi_{Y/S}$ we shall use the critical set
$C_\pi$ and the closure $D_\pi$ of $\pi(C_\pi)$, rather than the weak
discriminant set $D^{disc}_\pi$ (in case $D_\pi^{disc} \neq \pi(C_\pi)$).

\begin{pf} (1): Since $T_{X/Y}=0$ and $\pi^{-1}$ is exact we
have a commutative diagram of $\pi^{-1}(\Oc_Y)$-modules:
  \begin{displaymath} \xymatrix{0 \ar[r]& T_{X/S} \ar[r]&
T_{X/S \to Y/S} \ar[r]& \Cc_{X/Y} \ar[r]& 0\\ 0 \ar[r]&
\pi^{-1}(T_{Y/S}^\pi)\ar[u] \ar[r]&
\pi^{-1}(T_{Y/S})\ar[u]^{\psi_0} \ar[r]&
\pi^{-1}(\Ec_{X/Y})\ar[u]^{\phi_0}\ar[r]& 0 }
  \end{displaymath} Clearly, if $\psi_0$ is injective, then
$\phi_0$ is injective.

(2): Letting $\overline T_{X/S}$ be the image of $d\pi$ we have a
corresponding commutative diagram of $\Oc_Y$-modules
  \begin{displaymath} \xymatrix{0 \ar[r]& \pi_*
(\overline T_{X/S}) \ar[r]& \pi_*(T_{X/S \to Y/S}) \ar[r]&
\pi_*(\Cc_{X/Y}) \ar[r]& 0\\ 0 \ar[r]& T^\pi_{Y/S}\ar[u]
\ar[r]& T_{Y/S}\ar[u]^{\pi_*(\psi_0)} \ar[r]&
\Ec_{X/Y}\ar[u]^{\phi}\ar[r]& 0 }
  \end{displaymath} Here the map from $T_{Y/S}$ to
$\overline T_{X/S}$ is the composition $T_{Y/S}\to
\pi_*(T_{X/S})\to \pi_*(\overline T_{X/S})$ where by assumption
the latter map is surjective. By the definition of
$T^\pi_{Y/S}$ the assertion now follows.

  (3): The assumption implies $T_{X/S}= \overline T_{X/S}$, and
therefore the map $\phi$ is injective; that $\psi$ is
surjective implies that $\phi$ is surjective; therefore
$\phi$ is an isomorphism.

  The remaining assertion is evident, since if $\phi$ or
$\phi_0$ are injective, it follows that $\supp
\Ec_{X/Y}\subset \supp \Cc_{X/Y}$.
\end{pf}

Let $k\to A \xrightarrow{\pi} R$ be homomorphisms of rings,
and $T_{A/k}= Hom_A(\Omega_{A/k},A)$ and $T_{R/k}$ be the
Lie algebroid of $k$-linear derivations of $A$ and $R$,
respectively.  The tangent morphism $d \pi :T_{R/k} \to
T_{A/k\to R/k}= \Hom_A(\Omega_{A/k}, R)= \Hom_R(R\otimes_A
\Omega_{A/k} , R) $ is given by $d\pi (\partial)(da)=
\partial (\pi(a))$, $\partial \in T_{R/k}$ $a\in A$.  Let
$\psi_0 : T_{A/k} \to T_{A/k\to R/k}$ be the natural
morphism, which is injective if $\pi$ is injective. The
liftable derivations are $T_{A/k}^\pi = \Imo (d\pi)\cap
\psi(T_{A/k})$.
\begin{example}\label{ex:psi1} We continue with \Example{\ref{ex:psi}, (i)}.  We have $T_{A/k\to R/k}= R\frac 1t
  \partial_t$ and $\Cc_{R/A}= T_{A/k\to R/k}/T_{R/k}=
R/\mf_R =k$, but $T^\pi_{A/k}= T_{A/k}$, so $\Ec_{R/A}=0$.
\end{example}
Liftable derivations may be included in $\mf_R^NT_R$ for any high $N$.
We give a random example that $N$ can be high (in this case $N=9$).
\begin{example} Let $A$ be a $\Qb$-algebra with normalisation $R=
  \Qb[[t]]$, $A\neq R$.  Then $T_R^{cont} = T_{R/\Qb}^{cont}= R\partial_t$
  (continuous derivations; alternatively, work instead with the localisation
  $R=\Qb[t]_0$).  There exists a smallest integer $c$ such that $(t^c)\in
  A$. Put $R_c= R/(t^c)$ and $A_c= A/(t^c)$, so $A_c \subset R_c$ is an
  inclusion of finite-dimensional $\Qb$-algebras.  A derivation of $R$
  that preserves $A$ is of the form $\nabla = ft\partial_t$, $f\in R$, so $\nabla $
  also induces a derivation of $A_c$ and $R_c$. If $\bar a_1,\dots ,
  \bar a_l$ are generators of $A_c$, where $\bar a_i \in R_c$ is the
  projection of $a_i\in R$, the condition on $\bar f$ is $\bar f\bar
  a'_i \in A_c$, $i=1,2,3$, where $\bar a_i'= \overline{t\partial_t(a_i)}$.
  These equations may have only the trivial solution $f\in (t^c)$. A
  concrete example is $A= \Qb[[a_1, a_2, a_3, t^9, t^{10}, \dots]]$,
  $a_1=t^4+t^7, a_2=t^5+t^8, a_3=t^6$. We indicate that $T_A^\pi \subset
  \mf_R^9T_R$.  The ramification index is $4$ and $c=9$, so working in
  $R_{9}$ we see that $\bar 1, \bar a_1 , \bar a_2, \bar a_3$ forms a
  $\Qb$-basis for the $4$-dimensional subspace $A_{9}$ ($\bar a_i \cdot
  \bar a_j =0$, $i\neq j$) of the $9$-dimensional space $R_9$.  To solve
  $\bar f \bar a'_i \in A_9$ it suffices to make the ansatz $\bar f =
  c_0 + c_1t +c_2 t^2 + c_3t^3 + c_4 t^4\in R_9$, $c_i\in \Qb$, so we
  have the equations $(c_0 + c_1t +c_2 t^2 + c_3t^3 + c_4 t^4)a'_i \equiv
  \alpha_{i1}a_1 + \alpha_{i2} a_2 + \alpha_{i3} a_3 \mod ( t^9)$, $i=1,2,3$.
  Identifying coordinates in a basis of the $9$-dimensional space
  $R_9$ gives $5$ equations for each $i$, so there are $15$ equations.
  There are $5$ parameters in $\bar f$ and the $\alpha_{ij}$ give another
  $9$ parameters, giving $14$ parameters.  A straightforward
  computation shows that the only solution is $\bar f=0$.  Note that
  the vectors $(c_0 + c_1t + \cdots + c_4t^4)a_i'$ span a $5$-dimensional
  subspace while the vectors $\alpha_{i1}a_1 + \alpha_{i2} a_2 + \alpha_{i3} a_3$
  span a $3$-dimensional subspace of $R_9$ so we should expect that
  their intersection is $0$.
\end{example}
\subsection{Logarithmic derivations}\label{log-der} Let $k \to A
\xrightarrow{\pi} R$ be inclusions of rings, where we regard $A$ and
$R$ as algebras over $k$, e.g.  $k= \Zb$, the ring of integers, and
$R$ is a discrete valuation ring.  Assume that the extension of
fraction fields $K(R)/K(A)$ is separably algebraic, and therefore that
the $A$-module of derivations $T_{A/k}$ can be regarded as an $A$-Lie
sub-algebroid of $T_{K(R)/k}$, $T_{A/k} \subset T_{K(R)/k}$; put
$T^\pi_{A/k}= T_{A/k} \cap T_{R/k}$.

\begin{definition} Denote by $T_{A/k}(\log^\pi \mf_A)$ the
  $A$-submodule of $T_{A/k}$ of derivations $\partial$ such that $\partial (\phi )\in
  R \phi$ when $\phi \in \mf_A$, and denote by $T_{A/k}(\mf_A)$ the
  $A$-module of derivations such that $\partial (\mf_A)\subset \mf_A$.
\end{definition} 
The Lie algebroid $T_{A/k}(\log^\pi \mf_A)$ ($=T_{A/k}(\log^\pi \{x_1,
\dots  , x_r\})$, for a basis $\{x_1, \dots $ $, x_r\}$ of $\mf_A$) is a
sub-Lie algebroid of the Lie algebroid $T_{A/k}(\mf_A)$; elements in
the former module are {\it $\pi$-logarithmic}\/ and elements in the
latter are {\it logarithmic derivations}.

\begin{lemma}\label{log-der-lemma}Let $L\subset K(A)$ be a {\it subset}. Denote by
  $T_{A/k}(\log^\pi L)$ the $A$-submodule of $T_{A/k}$ of derivations
  $\partial$ such that $\partial (\phi )\in R \phi$ when $\phi \in L$.
\begin{displaymath}
  T_{A/k}(\log^\pi L)\subset  T_{A/k}(\log^\pi \mf_A) \subset T_{A/k}(\mf_A)
\end{displaymath}
where the latter is an inclusion of Lie algebroids.  If $L$ is either
a non-zero $A$-submodule of $K(A)$ or contains a basis of $\mf_A$,
then $T_{A/k}(\log^\pi L)= T_{A/k}(\log^\pi \mf_A)$.
\end{lemma} \underline{Note}: If $P$ is an $\mf_A$-primary ideal the
Lie algebroids $T_{A/k}(P)$ and $T_{A/k}(\mf_A)$ are in general not
included in one another, but if $\Char A=0$, then $T_{A/k}(P)\subset
T_{A/k}(\mf_A)$ (see e.g. \cite{kallstrom:preserve}). Note also that
for general subsets $L\subset K(A)$ we have that $T_{A/k}(L)=\{\partial \in T_{A/k} \
\vert \ \partial(\phi)\in L \text{ when } \phi\in L\}$ is a Lie
subalgebroid of $T_{A/k}$, but the $A$-module $T_{A/k}(\log^\pi L)$
need not be a Lie algebra.  Moreover, in general $T_{A/k}(\log^\pi L)
\not\subset T_{A/k}(L)$, also when $L$ is an ideal of $A$.

The global version of the above definition is as follows. Let $\pi: X\to
Y$ be a dominant generically separably algebraic morphism of schemes,
where $X$ is a Krull scheme, i.e. locally the spectrum of a Krull
ring. A Krull ring is an integral domain formed as an intersection $A=
\cap_\lambda R_\lambda$ of discrete valuation rings $R_\lambda$ in its fraction field
$K(A)$ such that every non-zero element in $K(A)$ is invertible in all
but finitely many $R_\lambda$\footnote{We will need only the first
  condition for our main results.}; a locally N{\oe}therian scheme is
Krull if and only if it is integral and integrally closed in its
fraction field.  Polynomial rings in infinitely many variables over a
field are non-N{\oe}therian Krull rings.  See
\citelist{\cite{bourbaki:commutative}*{Ch VII,\S1.3}
  \cite{nagata:local}*{Ch V} \cite{fossum:krull}} for a treatment of
Krull rings.

For a closed subset $D\subset Y$ the sheaf of {\it logarithmic
  derivations along $D$}\/ is $T_{Y/S}(I_D)=T_{Y/S}(\log D)= \{\partial \in T_{Y/S}\
\vert \
\partial (I_D)\subset I_D \}$, where $I_D$ is the (reduced) ideal of
$D$. Let $C= \cup\ C_i$ be a union of irreducible closed subsets of
pure codimension $1$ in $X$.  Let $\Oc_{X,c_i}$ be the local ring at
the generic point $c_i$ of $C_i$, and $\nu_i$ the associated
normalised discrete valuation.  The pre-sheaf of {\it $\pi$-logarithmic
  derivations}\/ along $C$ is
\begin{multline*} T_{Y/S}(\log^\pi C)= \{\partial \in
T_{Y/S} \ \vert \ \nu_{c_i}(\pi_{c_i}(\partial(f)))\geq
\nu_{c_i}(\pi_{c_i}(f))\\ \text{ when } f \in
\Oc_{Y,\pi(c_i)}, i=1, \dots , r \},
\end{multline*} 
where $\pi_c$ denotes the homomorphism of local rings $\Oc_{Y, \pi(c)}\to
\Oc_{X,c}$.  It is easy to see that $T_{Y/S}(\log^\pi C)$ actually is a
sheaf, and to see that a section $\partial$ of $T_{Y/S}$ is $\pi$-logarithmic
it suffices to check that its germ $\partial_{\pi(c_i)}\in
T_{\Oc_{Y,\pi(c_i)}/\Oc_{S,s_i}}(\log^\pi \mf_{\pi(c_i)})$, where $s_i\in
S$ is the image of $c_i\in Y$.  Put $D=\pi(C)$ and let $I_D$ be {\it
  any}\/ defining ideal of the closure of $D$, and put
$T_{Y/S}(\log^\pi I_D) = \{\partial \in T_{Y/S}\ \vert \ \nu_{c_i}(\partial (f))\geq
\nu_i(f)\quad \text{for all } \phi\in I, i=1, \dots ,r\}$. The notation is
incomplete since this sheaf also depends on the choice of $C$
(different $C$ can give the same $D$) but we choose not to burden the
notation; if $C= \pi^{-1}(D)$ there is no ambiguity and in practice $C$
will always be the critical set $C_\pi$ of the morphism $\pi$.  To
determine if $\partial\in T_{Y/S}$ belongs to $T_{Y/S}(\log^\pi I_D)$ it
suffices to check the condition for a set of local generators of
$I_D$.
\begin{lemma}\label{log-der-lemma2}
  \begin{enumerate}[label=(\theenumi)]
  \item The sheaf $T_{Y/S}(\log^\pi I_D)$ is independent of choice of
    defining ideal $I_D$ of the closure of $D$; more precisely
    \begin{displaymath}
      T_{Y/S}(\log^\pi I_D)
    = T_{Y/S}(\log^\pi C).
  \end{displaymath}
\item If $I_D$ is a defining ideal of $D$, then
    \begin{displaymath}
      T_{Y/S}(\log^\pi C) \subset T_{Y/S}(\sqrt I_{D});
    \end{displaymath}
    in particular $ T_{Y/S}(\log^\pi C) \subset T_{Y/S}(\log D)$.
\end{enumerate}
\end{lemma}
In view of $(1)$ it is reasonable to write
\begin{displaymath}
  T_{Y/S} (\log^\pi D) = T_{Y/S}(\log^\pi I_D) = T_{Y/S}(\log^\pi C),
\end{displaymath}
where $I_D$ is an
arbitrary defining ideal of the closure of $D$.
\begin{pf}
  (1): Clearly, $ T_{Y/S}(\log^\pi C)\subset T_{Y/S}(\log^\pi I_D)$,
  so assuming $\partial$ is a section of $ T_{Y/S}(\log^\pi I_D)$,
  regular at $\pi(c_i)$, and $f\in \Oc_{Y,\pi(c_i)}$, we have to check
  that $\nu_{c_i}(\partial(f))\geq \nu_{c_i}(f)$. If $f\notin
  \mf_{\pi(c_i)}$ we have $\nu_{c_i}(f)=0$ and the assertion is
  obvious. If $f\in \mf_{\pi(c_i)}$, then $f^n\in I_{\pi(c_i)}$ for
  sufficiently high $n$, hence $\nu_{c_i}(\partial(f^n)) \geq
  \nu_{c_i}(f^n)$, i.e.  $\nu_{c_i}(n\partial(f)f^{n-1}) =
  \nu_{c_i}(\partial (f)) + (n-1)\nu_{c_i}(\partial(f)) =
  n\nu_{c_i}(\partial(f))\geq n \nu_{c_i}(f)$, so $
  \nu_{c_i}(\partial(f))\geq \nu_{c_i}(f)$.

  (2): In the light of $(1)$ it suffices to prove that a section
  $\partial$ of $ T_{Y/S}(\log^\pi I_D)$ preserves the minimal
  associated primes of $I_D$ for any defining ideal of $D$.  A minimal
  associated prime of $I_D$ is of the form $\pi(c_i)$ for some generic
  point $c_i$ of $C$.  Therefore $\partial \in
  T_{Y/S,\pi(c_i)}(\log^\pi \mf_{\pi(c_i)})\subset
  T_{Y/S,\pi(c_i)}(\mf_{\pi(c_i)})$ \Lem{\ref{log-der-lemma}}.
\end{pf}


\subsection{Quasi-finite and algebraic morphisms}\label{quasi-finite-section}

\begin{prop} \label{zariski-main} Let $\pi: A \to R$ be an injective
  homomorphism of local integral rings, where $A$ is N{\oe}therian.
  Assume that at least one of the following conditions hold:
  \begin{enumerate}[label=(\roman*)]
  \item $A$ is complete.
  \item the $A$-algebra $R$ is of finite type.
  \item $R$ is integrally closed in $K(R)$, $K(R)$ is finite
over $K(A)$, and the integral closure of $A$ in $K(A)$ is
finite.
  \end{enumerate} Then the following are equivalent:
    \begin{enumerate}[label=(\theenumi)]
    \item $\dim_{k_A} k_A\otimes_A R < \infty $ ($R$ is
quasi-finite over $A$).
  \item $R$ is finite over $A$.
  \end{enumerate}
\end{prop} See \cite{matsumura} for a discussion of the problem with
the property that the normalisation be finite over $A$. It suffices
that the completion of $A$ be reduced, and if $R$ is a discrete
valuation ring this is also necessary.

\begin{remark}\label{ram-rem} Recall the situation when $A$
is a discrete valuation ring and $K(R)/K(A)$ is finite. Then
$k_R/k_A$ is finite, and the ramification index of $A\to R$
is the integer $e$ such that $R\mf_A = \mf_R^e$; it is also
the index of the value group of $A$ in the value group of
$R$.  If $A$ moreover is complete, then $[K(R):K(A)] =
e[k_R:k_A]$.
\end{remark}
\begin{pf} We have only to prove $(1)\Rightarrow (2)$. If (i) holds,
  see \cite[Ch. 0, Cor. 7.4.3]{EGA1}. If (ii) holds this is Zariski's
  main theorem, see \cite{EGA3}*{Ch. III, Cor. 4.4.6}. If (iii) holds
  it follows that the integral closure $A^*$ of $A$ in $K(R)$ is
  finite over $A$ \cite{bourbaki:commutative}*{Ch V,\S1.6, Prop.  18}
  and since $R$ is integrally closed, $A^* \subset R$. Now $R$ is a
  union of local quasi-finite $A^*$-algebras $F$ of finite type;
  $K(A^*)=K(R)$ since $K(R)/K(A)$ is algebraic; therefore
  $K(A^*)=K(F)$, so by Zariski's main theorem $F=A^*$ \cite{EGA3}*{Ch. III, Cor. 4.4.8}; hence $R= A^*$.
\end{pf}

\begin{prop}\label{q-finite} Let $\pi: A \to R$ be a local
injective homomorphism of N{\oe}therian integral rings.
Consider the conditions:
    \begin{enumerate}[label=(\theenumi)]
  \item $K(R)/K(A)$ is algebraic.
  \item $k_R/k_A$ is algebraic.
  \end{enumerate} Then:
  \begin{enumerate}[label=(\alph*)]
  \item If $(1)$ and $(2)$ holds, and $A$ is universally
catenary (e.g.  formally equidimensional), then $\dim A =
\dim R$.
\item If $R$ is a valuation ring and $A= K(A)\cap R$, then
$(1)\Rightarrow (2)$.
\item \label{finite-q-finite} If $F\mf_A$ is $\mf_F$-primary
for each local intermediate $A$-algebra $(F,\mf_F)$ of
finite type, $A\subset F \subset R$, then $(2)\Rightarrow
(1)$.
  \end{enumerate}

\end{prop} Assuming (2) the assumption in (c) is the same as
$F/A$ is quasi-finite.  Note that if $R/A$ is of finite type
and $F$ is a local $A$-algebra as in~\ref{finite-q-finite},
then if $R/A$ is quasi-finite it follows that $F/A$ is
quasi-finite. To see this, apply $k_A\otimes_A \cdot $ to
the exact sequence of $A$-modules $0 \to F \to R \to F/R\to
0$, noting that $\dim_{k_A} Tor_1^A(k_A , R/F) <\infty$ when
$R$ is of finite type, and by assumption
$\dim_{k_A}k_A\otimes_A R< \infty $. If $R/A$ is not of
finite type this does not follow (example:
$A=k[x]_{0}\subset R= k[[x]]$).

\begin{pf} (a) An algebraic field extension $K(R)/K(A)$ is a
union of algebraic extensions $K(F)/K(A)$ where $A\subset
F\subset R$ and $F$ is of finite type over $A$.  Since $R$
is N{\oe}therian $\dim R = \sup_{A\subset F\subset R} \dim
F$; hence it suffices to see that $\dim F = \dim A$. Since
$F/A$ is of finite type we conclude from Ratliff's dimension
equality (see \cite{matsumura}*{Th.  15.6}), since $A$
is N{\oe}therian and universally catenary, that $(1-2)$
implies $\dim A = \dim R$.

  (b) See \cite{bourbaki:commutative}*{Ch VI, \S8, Prop. 1}.

  (c) By assumption $R$ is a union of local $A$-algebras $F$
of finite type such that $F\mf_A$ is $\mf_{F}$-primary,
hence by $(2)$ $F\otimes_A k_A$ is finite over $k_A$; hence
by Zariski's main theorem \Prop{\ref{zariski-main}} $F$ is
finite over $A$.  Since $K(R)$ is a union of the fields
$K(F)$ it follows that $K(R)$ is algebraic over $K(A)$.
\end{pf}

\subsection{Ramification}\label{ram-section}

We collect some facts about ramifications of injective
homomorphisms of local rings $k\to A \to R $. The cotangent
mapping is denoted
\begin{displaymath} \phi: k_R \otimes_{k_A}\mf_A/\mf_A^2 \to
\mf_R/\mf_R^2.
\end{displaymath}
A field extension $k\subset l$ has a finite differential
basis if $\dim_l \Omega_{l/k}<\infty$, which is equivalent
to $\trdeg l/k <\infty $ when $\Char k =0$ and equivalent
to $l/k$ having a finite $p$-basis when $\Char k =p$ (see
\cite[\S{} 26]{matsumura}). We say that an extension $k\subset
l$, where $k$ is a subring of the field $l$, has a finite
differential basis if $l$ has a finite differential basis
over the fraction field of $k$.
\begin{prop}\label{basic-prop} Let $(A, \mf_A, k_A)$ be a
local N{\oe}therian subring of a discrete valuation ring
$(R,\mf_R, k_R)$, where $\mf_A \neq 0$.  Let $k \subset A$
be a sub-ring such that the extensions $k_R/k$ and $k_A/k$
are separable with finite differential bases.  Consider the
following conditions:
    \begin{enumerate}[label=(\theenumi)]
  \item The map $\phi$ is injective.
  \item The ring $(A, \mf_A, k_A)$ is a discrete valuation
ring and $T_{A/k}= T^\pi_{A/k}$.
  \item The tangent morphism $d\pi: T_{R/k} \to T_{A/k\to
R/k}$ is surjective.
  \item The map $\phi$ is surjective.
  \end{enumerate} Then $(1)\Leftrightarrow (2)$,
$(1-2)\Rightarrow (3)$ and $(1-2)\Rightarrow (4)$. If $A/k$
is smooth, then $(3)\Rightarrow (1-2)$. If $A$ is a discrete
valuation ring, then $(4)\Rightarrow (1-2)$.
\end{prop}
\begin{pf} $(1)\Leftrightarrow (2)$: Since $A$ is a
N{\oe}therian integral domain, we have by (1) $1 \leq \dim A
\leq \dim_{k_A} \mf_A/ \mf_A^2$. Hence (1) implies that $(A,
\mf_A, k_A)$ is a discrete valuation ring, thus in either
case we have an inclusion of discrete valuation rings
$A\subset R$.  Since $\dim_{k_R}\Omega_{k_R/k}< \infty $ and
$\dim_{k_A}\Omega_{k_A/k}< \infty$ it follows from the
second fundamental exact sequence of differentials that
$\Omega_{R/k}$ and $\Omega_{A/k}$ are free of finite rank.
Therefore $T_{R/k}$ and $T_{A/k \to R/k} =
R\otimes_AT_{A/k}$ (see proof of \Proposition{\ref{psi}})
are free of finite rank over $R$.  It follows by Nakayama's
lemma that the map $d\pi : T_{R/k}\to R\otimes T_{A/k}$ is
surjective if and only if the tangent map $k_R \otimes_R
T_{R/k}\to k_R\otimes _{R}T_{A/k}$ is surjective, i.e.\  the
map $(\mf_R/\mf_R^2)^* \oplus T_{k_R/k} \to k_R\otimes_k
(\mf_A/\mf_A^2)^* \oplus k_R\otimes_{k_A}T_{k_A/k}$ is
surjective.  Since $k_R/k_A$ is separable, the map
$T_{k_R/k}\to T_{k_A/k \to k_R/k}=
k_R\otimes_{k_A}T_{k_A/k}$ is surjective; hence $d\pi$ is
surjective if and only if the map $(\mf_R/\mf_R^2)^* \to
k_R\otimes_{k_A}(\mf_A/\mf_A^2)^* = Hom_{k_A}(\mf_A/\mf_A^2,
k_R)$ is surjective and this map is surjective if and only
if its image contains the image of the natural map
$Hom_{k_A}(\mf_A/\mf_A, k_A )\to Hom_{k_A}(\mf_A/\mf_A, k_R
)$, i.e. $d\pi$ is surjective if and only if $T_{A/k}\subset
d\pi (T_{R/k})$.  Since (1) is equivalent to
$(\mf_R/\mf_R^2)^* \to k_R\otimes_{k_A}(\mf_A/\mf_A^2)^*$
being surjective, it follows that $(1)\Leftrightarrow (2)$.

  $(1-2)\Rightarrow (3)$: Since $\Omega_{A/k}$ is (free) of
finite type the map $d\pi$ is surjective if and
only if the map $ (\mf_R/\mf_R^2)^* \oplus T_{k_R/k} \to k_R
\otimes_{k_A} (\mf_A/\mf_A^2)^* \oplus
k_R\otimes_{k_A}T_{k_A/k}$ is surjective.  The assertion is
then clear from the above argument, as is the assertion
$(3)\Rightarrow (1-2)$ when $A/k$ is smooth for then
$\Omega_{A/k}$ is free of finite rank since
$\dim_{k_A}\Omega_{k_A/k}< \infty$.

  $(1-2)\Rightarrow (4)$: If $\phi$ is injective it is also
surjective since $\dim \mf_R/\mf_R^2 =1$. If $A$ and $R$ are
discrete valuation rings then $\phi$ is injective if and
only if it is surjective.
\end{pf}

\begin{prop}\label{etale-field} The following are equivalent
for a field extension $l/k$:
    \begin{enumerate}[label=(\theenumi)]
    \item $l/k$ is formally \'etale and algebraic.
  \item $l/k$ is separably algebraic.
  \item $\Omega_{f/k}=0$ for each intermediate field $l/f/k$
such that $f/k$ is of finite type.
  \end{enumerate}
\end{prop} If $l/k$ is of finite type {\it or} $l/k$ is separable,
then (3) is equivalent to $\Omega_{l/k}=0$.
\begin{example} Let $k$ be a field of characteristic $p>0$
and $l= k(x^{p^{-\infty} }) = k(x, x^{p^{-1}}, x^{p^{-2}},
\dots )$. Then $l/k$ is separable and formally \'etale, but
not algebraic (see \cite[Ex.  26.8]{matsumura}); moreover
$\Omega_{l/k} =0$, and there exist finite intermediate
fields $l/f/k$ such that $\Omega_{f/k}\neq 0$.
\end{example}

$(3)\Rightarrow (2)$ is proven in \cite{EGA4:1}*{Cor.
21.7.4} when $l/k$ is of finite type, and the argument below
is not very different, but it may be convenient to have it
presented in this way.
\begin{pf} $(1)\Leftrightarrow (2)$: see \cite{matsumura}*{Thms. 25.3,
    26.9}. For the remaining parts let $\Omega_k$ and $\Omega_l$
  denote differentials over the prime field and consider the exact
  sequence $0 \to \Gamma \to l\otimes\Omega_{k} \xrightarrow{h}
  \Omega_{l}\to \Omega_{l/k}\to 0$.  $(2)\Rightarrow (3)$: Then $l/k$
  is formally \'etale and $h$ is an isomorphism [{\it loc.\ cit}].
  $(3)\Rightarrow (2)$: If $k \subset f \subset l$ is an intermediate
  field such that $f/k$ is of finite type we have $\Omega_{f/k} = 0$
  and replacing $l$ by $f$ in the above exact sequence we see that the
  map $h$ is surjective and by Cartier's equality $\rko_f \Omega_{f/k}
  = \trdeg_k f + \rko_f \Gamma$ it follows that $\trdeg_k f =0$ and
  $\rko_f \Gamma =0$, hence $f/k$ is algebraic (hence finite) and $h$
  is injective, implying that $f/k$ is separable [loc.\ cit., Th.
  26.6].  Since $l$ is a union of such subfields $f$ it follows that
  $l/k$ is separably algebraic.
\end{pf}

\begin{thm}\label{formunramified} Let $(A,\mf_A, k_A)
\subset (R, \mf_R, k_R)$ be an inclusion of local
N{\oe}therian rings. Consider the conditions:
    \begin{enumerate}[label=(\theenumi)]
    \item $\Omega_{R/A}=0$.
  \item $k_R/k_A$ is separably algebraic and $R\mf_A =
\mf_R$.
  \item $k_R/k_A$ is separably algebraic and the map $\phi$
is surjective.
  \item $k_R/k_A$ is separably algebraic and the map $\phi$
is injective.
  \end{enumerate} Then $(2)\Leftrightarrow (3)$ and
$(2)\Rightarrow (1)$.  If the residue field extension
$k_R/k_A$ is finitely generated, then $(1)\Rightarrow (2)$.
If $A$ and $R$ are discrete valuation rings, then $
(2)\Leftrightarrow (3)\Leftrightarrow (4)$.
\end{thm}
Auslander and Buchsbaum \cite{auslander-buchsbaum} initiated the
systematic ramification theory of N{\oe}therian rings.  They proved
$(1)\Leftrightarrow (2)$ under the assumption that $J=\Ker (R\otimes_A
R\to R, r_1\otimes r_2 \mapsto r_1 r_2)$ is a finitely generated ideal
in $R\otimes_A R$, and hence that $\Omega_{R/A}$ is $R$-finite.  In
\cite{EGA4:4}*{\S17, Th. 17.4.1} one finds a different proof of the
same assertion assuming $R$ is a, possibly non-N{\oe}therian, finitely
presented $A$-algebra (implying $J$ is finitely generated), but the
argument in [{\it loc.\ cit.}]  for $(1)\Rightarrow (2)$ is fairly
involved.  In \Theorem{\ref{formunramified}} $J$ is not assumed to be
finitely generated.

\begin{pf} We give the proof only when $A$ is of equal characteristic.
  When $A$ is of mixed characteristic one replaces the residue fields
  by quasi-coefficient rings in the diagram below; these coefficient
  rings are formally smooth over the integers.  Thus we assume that
  $A$ and $R$ contain the prime field of $k_A$. Letting $\Omega_A$ and
  $\Omega_R$ denote differentials over the prime field we have the
  exact sequence $0 \to \Gamma\to R\otimes \Omega_{A}\xrightarrow{h}
  \Omega_{R}\to \Omega_{R/A}\to 0$.

  $(1)\Rightarrow (2)$: Since $k_R$ and $k_A$ are separable and hence
  formally smooth over the prime field \cite[Th. 26.9]{matsumura} we
  have split exact sequences $0 \to \mf_A/\mf_A^2 \to
  k_A\otimes_A\Omega_A \to \Omega_{k_A}\to 0$ and $0 \to \mf_R/\mf_R^2
  \to k_R\otimes_R\Omega_R \to \Omega_{k_R}\to 0$ and therefore a
  commutative diagram which is exact also to the left:
  \begin{displaymath} \xymatrix{0\ar[r]&
k_R\otimes_{k_A}\mf_A/\mf_A^2 \ar[r]\ar[d]^\phi&
k_R\otimes_A \Omega_A\ar[r]\ar[d]^f&
k_R\otimes_{k_A}\Omega_{k_A}\ar[d]^g\ar[r]&0\\ 0\ar[r]&
\mf_R/\mf_R^2\ar[r]& k_R\otimes_R \Omega_R\ar[r]&
\Omega_{k_R}\ar[r]&0.  }
\end{displaymath} By (1) the map $f$ is surjective, hence
$g$ is surjective, so $\Omega_{k_R/k_A}=0$.  By Cartier's
equality, using $k_R/k_A$ is finitely generated,
\begin{displaymath}
  \rko_{k_R}\Omega_{k_R/k_A}= \trdeg_{k_A}k_R +
\rko_{k_R}\Ker g ,
\end{displaymath}
it follows that $\Ker g =0$ and $\trdeg_{k_A}k_R =0$, implying that
$k_R/k_A$ is separably algebraic [{\it loc.\ cit.}, Th. 26.6].
Therefore $g$ is an isomorphism, in particular injective, hence, as
$f$ is surjective, $\phi$ is surjective by the ``serpent lemma''; by
Nakayama's lemma it follows that $R\mf_A = \mf_R$.

  $(2)\Rightarrow (1)$: Assuming (2) the ring $R$ is a union
of local $A$-algebras $F$ of finite type, $A \subset F
\subset R$, such that $k_F/k_A$ is separably algebraic and
$F\mf_A = \mf_F$.  It is easy to see that (1) will follow if
$\Omega_{F/A}=0$ for all such $F$.  Replacing $R$ by $F$ in
the above commutative diagram, by (2) the map $g$ is an
isomorphism and $\phi$ is surjective, hence $f$ is
surjective, hence, since $\Omega_{F/A}$ is of finite type,
by Nakayama's lemma $\Omega_{F/A}=0$.

  $(2) \Leftrightarrow (3)$ is evident.

  If $A$ and $R$ are discrete valuation rings, then $\phi$
is an isomorphism if and only if it is injective or
surjective, implying $(3)\Leftrightarrow (4)$.
\end{pf}

The study of a ramification of a morphism $X\to Y$ can often be
reduced to the case when $X$ is the spectrum of a discrete valuation
ring. We now decide what a ramified morphism is in this case.

\begin{definition} \label{tame-ramification} 

  Let $\pi: A\to R$ be an injective homomorphism of local
  N{\oe}therian rings of the type in \Proposition{\ref{basic-prop}}.
  Then $\pi$ is {\it ramified}\/ if the two equivalent conditions
  (1-2) in \Proposition{\ref{basic-prop}} do not hold.  A {\it tamely
    ramified}\/ morphism $\pi$ is a ramified morphism such that the
  induced fraction field extension $K(R)/K(A)$ is separably algebraic,
  the residue field extension $k_R/k_A $ is separable, and if $x\in
  \mf_A \setminus \mf_A^2$, then $\nu(\pi(x))$ is co-prime to the
characteristic of $A$, where $\nu$ is the normalised valuation of $R$.
A local domain $A$ of unequal characteristic with $\Char k_A =p>0$, is
(absolutely) unramified if $p\notin \mf_A^2$, and tamely ramified if
$p\in \mf_A^n\setminus \mf_A^{n+1}$ for some integer $n\geq 2$ such
that $p{\not|}\ n$.  A morphism of schemes $\pi: X \to Y$, where $X$
is regular in codimension $\leq 1$, is tamely ramified if for each
point $x$ of height $1$ in $X$ the map of local rings
$\Oc_{Y,\pi(x)}\to \Oc_{X,x} $ is tamely ramified. An integral scheme
$Y$ of unequal characteristic is unramified (tamely ramified) if each
local ring $\Oc_{Y,y}$ is unramified (tamely ramified).
\end{definition} When $\Char k_A=0$ ramified morphism are
always tamely ramified. If $R/A$ is not ramified, then we
say that $R/A$ is {\it unramified}, and if a ramified
morphism is not tame, then it is {\it wild}.

\begin{remark}
  \begin{enumerate}[label=(\roman*)]
  \item Assuming $K(R)/K(A)$ and $k_R/k_A$ are separably
finite, our definition of tameness is not the same as
Abhyankar's \cite[\S10, p. 245]{abhyankar:res-surf} when
$A$ is of mixed characteristic.
\item \label{form-unr} Recall that an $A$-algebra $R$ is formally
  unramified if and only if $\Omega_{R/A}=0$. This is not the same as
  unramified in \Definition{\ref{tame-ramification}}, essentially
  because $\phi$ in \Proposition{\ref{basic-prop}} need not be
  injective when it is surjective.  Assuming $R$ is a discrete
  valuation ring, by \Theorem{\ref{formunramified}} and
  \Proposition{\ref{basic-prop}} it follows that if $R/A$ is
  unramified, then $R/A$ is formally unramified; conversely, if $R/A$
  is a formally unramified extension of discrete valuation rings, and
  $k_R/k_A$ is finitely generated, then $R/A$ is unramified.
  \end{enumerate}
\end{remark}
If $R/A$ has wild ramification we may have $T^\pi_{A/k}\subset \mf^NT_{A/k}$
for high $N$.
\begin{example} Let $k$ be a field of characteristic $p$ and
$\pi: A\cong R =k[y]\to R=k[x]$ be a morphism over $k$. If
$y=x^n$, and $n$ is co-prime to $p$, then $T_{A/k}^\pi =
Ay\frac{d}{dy}$. If $y= x^{p}/(1-x^{p-1})$, then
$T^\pi_{A/k} = Ay^2\frac{d}{dy}$. One can iterate the latter
example, letting $\pi$ be the composition $A\cong R \to R\to
\cdots \to R$ ($r$ times). This gives $T^\pi_{A/k} = A
y^{2r} \frac{d}{dy}$.
\end{example}

\begin{lemma}\label{p-rings} Let $\phi_0 :k\to k'$ be a separably
  algebraic extension of fields.
  \begin{enumerate}[label=(\theenumi)]
  \item A homomorphism $\psi: k[[X]]\to k'[[X']]$ inducing $\phi_0$ on
    residue fields is $(X')$-\'etale if and only if it is
    $(X')$-unramified if and only if $\psi (X)\in (X')\setminus
    ({X'}^2)$. \\
  \item If $\Char k= p >0$ then a local homomorphism of complete
  $p$-rings $\phi : (W(k), tW(k), k) \to (W(k'), t'W(k'), k')$,
  extending $\phi_0$, is $(t')$-\'etale if and only it is
  $(t')$-unramified if and only if $\phi (t)\in (t')\setminus (t'^2)$.
  \end{enumerate}
\end{lemma} 

\begin{pf} Let $R(k)$ denote $k[[t]]$ if $\Char k =0$ and if $\Char k
  =p >0$, it denotes $k[[t]]$ or the $p$-ring $W(k)$.  We have a
  homomorphism $\phi:(R(k), tR(k), k)\to (R(k), sR(k), k)$, inducing
  $\phi_0 : k\to k'$ on residue fields. If $\phi(t)\in (s^2)$ the
  morphisms is not unramified in the $(s)$-adic topology, so it
  suffices to see that if $\phi(t)\in (s)\setminus (s^2)$, then $\phi$
  is \'etale in the $(s)$-adic topology.  Since it is clearly
  $(s)$-unramified it remains to see that $\phi$ is $(s)$-smooth.
  Since in either case $R(k')$ is an integral domain containing the
  valuation ring $R(k)$ it is flat over $R(k)$; the fibre
  $R(k')/tR(k')= k'$ is separably algebraic over $k$, hence formally
  \'etale; therefore by \cite{EGA4:1}*{Th.  19.7.1} $R(k')$ is
  $(s)$-\'etale over $R(k)$.
\end{pf}

\begin{remark} Assume that $k$ is perfect so we have unique inclusions
  $W(f) \subset W(g)$ over field extensions $k'/g/f/k$, and therefore
  an unambiguous definition of $ \dirlim_{k'/f/k} W(f)$; this ring is
  formally \'etale over $W(k)$ (see argument e.g. in
  \cite{spivakovsky:popescu}).  However, the canonical injective
  homomorphism $\dirlim_{k'/f/k} W(f)\to W(\dirlim_{k'/f/k}f)= W(k')$,
  where the limits are over finite intermediate fields $f$, is not
  surjective when $k'/k$ is not finite, and there exist no
  intermediate ring $\dirlim_{k'/f/k} W(f) \subsetneq R \subset W(k')$
  such that $\dirlim_{k'/f/k}W(f)\to R$ is finite \'etale; therefore
  $W(k')/W(k)$ is not an inductive limit of finite \'etale morphisms.
  By N\'eron desingularisation \cite{neron:minimaux} (recall also
  Popescu's generalisation in \cite{popescu:neron} and its new proof
  in \cite{spivakovsky:popescu}), formally smooth morphisms of discrete
  valuation rings are filtered inductive limits of finite type smooth
  morphisms; hence formally \'etale morphisms are inductive limits of
  (finite) \'etale morphisms; therefore $W(k')$ is {\it not}\/
  formally \'etale over $W(k)$ when $k'/k$ is not finite.  Actually,
  letting $R(k)$ be as in the proof and the morphism $R(k)\to R(k')$
  be $(s)$-unramified, then it is formally \'etale if and only if it
  is finite \'etale, and this happens if and only if $k'/k$ is
  separably finite.
\end{remark}

\section{Liftable derivations of separably algebraic
morphisms}\label{GF}

\subsection{Lifting derivations to discrete valuation rings} This
section contains the main step in the description of liftable
derivations.  Let $\pi :(A, \mf_A, k_A) \to (R, \mf_R, k_R)$ be an
injective homomorphism of local N{\oe}therian algebras over a ring $k$,
where $R$ is a discrete valuation ring, and assume that $A$ and $\pi$
are tamely ramified \Defn{\ref{tame-ramification}}.  In particular,
derivations of $A$ are uniquely liftable to derivations of the
fraction field of $R$.
\begin{thm}\label{fund-thm} Let $\pi: A\to R$ be a local injective
  homomorphism of N{\oe}therian rings, where $R$ is a discrete
  valuation ring and the extension of fraction fields $K(R)/K(A)$ is
  separably algebraic.
    \begin{enumerate}[label=(\theenumi)]
    \item If $\pi $ is tamely ramified, then
    \begin{displaymath} T^\pi_{A/k} \subseteq
T_{A/k}(\log^\pi K(A))\cap T_{R/k}(\mf_R)\cap T_{A/k}(\mf_A).
    \end{displaymath}
  \item Assume that $\pi$ is residually algebraic and tamely ramified,
    and if $A$ is of unequal characteristic then $A$ is tamely
    ramified. Then
    \begin{displaymath} T^\pi_{A/k}= T_{A/k}(\mf_A).
    \end{displaymath}
  \item If $K(A)= K(R)$, then
    \begin{displaymath} T_{A/k}(\log^\pi \mf_A)\subset
T_{A/k}(\mf_A)\cap T_{R/k}(\mf_R).
\end{displaymath}
  \item  \label{main-part} If $\pi$ is tamely ramified and  $\mf_A \subseteq  L\subseteq  K(A)$,
then
\begin{displaymath} T_{A/k}(\log^\pi L)= T^\pi_{A/k}\cap T_{A/k}(\mf_A)=
  T_{A/k}^\pi.
\end{displaymath} If moreover $\pi$ is residually
algebraic and $A$ is tamely ramified (when $A$ is of unequal
characteristic), then
    \begin{displaymath} T^\pi_{A/k} = T_{A/k}(\mf_A) =
T_{A/k}(\log^\pi \mf_A).
    \end{displaymath}
  \end{enumerate}
\end{thm}

The main part of \Theorem{\ref{fund-thm}} is ~\ref{main-part}. There
is a certain stability built in in this result, so that the liftable
derivations do not change if $R$ is modified by a tamely ramified
morphism. Let $\pi: A \xrightarrow{\pi_1} R \xrightarrow{\pi_2} B$, where
$B$ also is a discrete valuation ring.  Assuming that $\pi_1$ is
residually algebraic and $\pi_2$ is tamely ramified we have
\begin{displaymath} T_{A/k}^\pi = T_{A/k}(\log^{\pi_1} \mf_A) \subset
  T_{R/k}(\mf_R) = T_{R/k}^{\pi_2} \subset T_{B/k},
\end{displaymath} and therefore  $T_{A/k}^{\pi} =
T_{A/k}^{\pi_1}$.
\begin{lemma}\label{fund-lemma} Let $(A, \mf_A, k_A)$ be a
local N{\oe}therian domain and $\partial$ a derivation of
$A$. Let $\nu$ be an additive valuation on the fraction
field $K(A)$ with a centre in $A$.  There exists an element
$c$ in the (ordered) value group of $\nu$ such that for
$x\in K(A)$
  \begin{displaymath} \nu(\frac{\partial(x)}x) \geq c.
  \end{displaymath}
\end{lemma}
\begin{pf} Select generators $\{x_1, \dots ,x_r\}$ of the
maximal ideal $\mf_A$.  Let $\partial$ also denote its
unique extension to a derivation of $K(A)$.  Let $c$ be an
element in the value group satisfying $\nu (\partial
(x_i)/x_i) \geq c$, $i=1, \dots , r$, where we may assume
$c\leq 0$.  Then if $x= ab \in K(A)$, where $a,b\in K(A)$,
one has $ \partial (x)/x = \partial (a)/a +\partial(b)/b$;
hence
  \begin{equation}\label{fund-rel}
\nu(\frac{\partial(x)}x)\geq \min \{\nu(\frac{\partial
(a)}a), \nu(\frac{\partial (b)}b)\}.
  \end{equation} If $a\in A$ we have $a=
x_{i_1}x_{i_2}\cdots x_{i_n}u$ where $x_{i_j}\in \{x_1, x_2,
\dots , x_r\}$, and $u\in A\setminus \mf_A$ is a unit. Then
$\nu(\partial (u)/u) \geq 0 $ and by (\ref{fund-rel}), $\nu
(\partial(a)/a)\geq c$. Again by (\ref{fund-rel}) we get for
any $x=a/b\in K(A)$ that $\nu(\partial(x)/x)\geq c$.
\end{pf}

\begin{pfof}{\Theorem{\ref{fund-thm}}} Put $p= \Char R$.  (1): We
  first prove $T_{A/k}^\pi \subset T_{R/k}(\mf_R)$, so letting
  $\partial$ be a derivation of $K(A)$, hence also a derivation of
  $K(R)$ since $\pi$ is generically separably algebraic, such that
  $\partial (R)\subset R$ and $\partial(A)\subset A$, we have to see
  that $\partial (\mf_R)\subset \mf_R$.
  
Assume the contrary, that $\partial (\mf_R)\nsubseteq \mf_R$. Let
$s\in \mf_R \setminus \mf_R^2$ be a uniformising parameter, so any
element $x\in R$ is of the form $x= us^n$ for a unit $u$ and positive
integer $n$.  From the relation
  \begin{equation}\label{m-realisation}
    \partial(x)=(\frac{\partial(u)}{u} + n\frac{
    \partial (s)}{s})x
\end{equation} and the assumption it follows  that $\partial
(s)$ must be  a unit. Let $\nu$ be the normalised valuation
of $R$, so $\nu (s^n)=n$, $\nu (\partial (s)) =0$.
Let $t \in \mf_A \setminus \mf_A^2$ have smallest possible value $n=
\nu (t)$. Since by tameness $p{\not \vert \ }n $ it follows that
$\partial (t)\neq 0$ and $\nu (\partial (t)) = n-1$, hence $\partial
(t)=u$ is a unit in $A$ and $n=1$; thus $t\in A$ is a uniformising
parameter of $R$.  Replacing $\partial $ by $u^{-1}\partial $ we can
assume that $\partial (t) =1$.  Let $\phi : R \to k_R$ be the
projection to the residue field of $R$.  We regard $k_A$ as a subfield
of $k_R$, using the natural inclusion map. Let $\{t, x_1, \dots ,
x_r\}$ be a subset of $\mf_A$ that induces a basis of $\mf_A/\mf_A^2$.
Since $R/A$ is ramified we have $r\geq 1$, and writing $x_i = u_i
t^{n_i}= u_it^{\nu(x_i)}$ , $u_i\in R\setminus \mf_R$, then $\phi (u_i)\notin k_A$.
Consider the set
  \begin{displaymath} S =\{x\in A \ \vert \ \phi (\frac
{x}{t^{\nu(x)}}) \in k_R \setminus k_A\},
\end{displaymath} which is nonempty since each $x_i$
belongs to $S$.  If $x= ut^n$, for some $u\in R \setminus
\mf_R$, and $p{\not \vert \ } n$, by (\ref{m-realisation}), we have $\nu (\partial (x))= n-1$
and hence $\phi (\partial (x)/t^{n-1}) = n \phi (u)$.  This
gives the following relation in $k_R$ for $x\in A$ such that
$p {\not\vert \ } \nu(x)$ (so in particular $\nu(x)\neq 0$)
  \begin{equation}\label{value} \phi (\frac{\partial
(x)}{t^{\nu (\partial (x))} })= \nu(x) \phi (\frac
{x}{t^{\nu(x)}}).
\end{equation} Put $n=\nu(S)$, the minimal value of the
elements in $S$. Then $n\geq 1$, for if $x\in A$ and
$\nu(x)=0$, then $\phi(x/t^{\nu(x)}) = \phi (x)\in k_A$, so
$x\notin S$.  Select $x\in S$ such that $\nu(x) =n $.  Then
clearly $x\in \mf_A \setminus \mf_A^2 $, hence $p {\not
  \vert \ }\nu(x)$ since $R/A$ is tamely ramified, hence $\nu (\partial (x)) = n-1$; by
(\ref{value}) we get that $
\partial (x)\in S$, contradicting the minimality of $n$.  This completes the proof  that
$\partial(\mf_R)\subset \mf_R$.

Since $T_{A/k}\cap T_{R/k}(\mf_R) = T_{A/k}(\mf_A)$, we get
$T_{A/k}^\pi \subset T_{A/k}(\mf_A)$.  It remains to see that
$T_{A/k}^\pi \subset T_{A/k}(\log^\pi L)$. Let $a\in L$, so $a= vt^m$
where $v$ is a unit in $R$, and let $\partial \in T_{A/k}^\pi$. We know
that $\partial (t)\in Rt$, hence by (\ref{m-realisation})
  \begin{displaymath} \frac {\partial (a)}a = \frac
{\partial (v)}v + m \frac {\partial(t)}t \in R.
\end{displaymath} Therefore $\partial \in T_{A/k}(\log^\pi
L)$.

\renewcommand{\theenumi}{\arabic{enumi}}\renewcommand{\labelenumi}{\theenumi.}
(2): Continue to regard, by generic separable algebraicity,
derivations $\partial \in T_{A/k}(\mf_A)$ as derivations of $K(R)$ and in
particular as derivations $R \to K(R)$. We need to prove that $\partial (R)\subset
R$.  Let $\nu$ be the normalised valuation of $K(R)$ whose valuation
ring is $R$.  By \Lemma{\ref{fund-lemma}} there exists an integer $c$
such that $\nu(\partial(r))\geq \nu(r) -c$ for $r\in R$, hence $\partial(\mf_R^n)\subset
\mf_R^{n-c} \subset K(R)$.  Therefore, providing $K(R)$ with the linear
topology where $\{\mf_R^n\}$ is family of neighbourhoods of $0$, it
follows that $\partial $ is a continuous derivation $R \to K(R)$; hence it
has a unique continuous extension to the completed rings $\partial :\bar R
\to \varprojlim_n K(R)/\mf_R^n \cong K(\bar R)$; hence it has a unique
continuous extension to a derivation $K(\bar R)\to K(\bar R)$.  By
Artin-Rees' lemma $R\subset \bar R$, where $\bar R$ also is a discrete
valuation ring \cite[ VI, \S{} 5.3, Prop 5 ]{bourbaki:commutative}.  Let
$r\in R$. If $\partial (r)\notin R$, then $1/\partial (r)\in \mf_R \subset \mf_{\bar R}$, so
$\partial (r)\notin \bar R$.  Therefore it suffices to see that $\partial (\bar R)\subset
\bar R$. The map $\pi: A\to R$ induces a map $\bar \pi :\bar A \to \bar R$
of completed rings. Since $\bar R$ is integral and $\bar A$ is
1-dimensional \Lem{\ref{q-finite}} $\Ker \bar \pi$ is either $\mf_{\bar
  A}$ or $0$; since $A\subset \bar A$ and $R\subset \bar R$ by Artin-Rees' lemma
it follows that $\bar \pi (A) \neq 0$ and therefore $\bar \pi$ is
injective (cf. \Remark{\ref{gabrielov}}). Since $T_A(\mf_A)\subset T_{\bar
  A}(\mf_{\bar A})$ it should be clear now that one may assume that
$A= \bar A$ and $R = \bar R$.

  There are two cases:

  (a) $A$ and $R$ are complete of unequal characteristic: Put $p=\Char
  k_A = \Char k_R $.  The rings $A$ and $R$ contain complete $p$-rings
  $W(k_A)\subset A$, $W(k_R)\subset R$, and the field extension
  $k_R/k_A$ lifts to a formally unramified inclusion of $p$-rings
  $(W(k_A), tW(k_A)\\ , k_A) \subset (W'(k_R), sW'(k_R), k_R)$ taking $t$
  to $s$ \cite{matsumura}*{Thms.  29.1-2}, which is hence
  $\mf_{W'(k)}$-\'etale \Lem{\ref{p-rings}}. Since all complete
  $p$-rings with given residue field are isomorphic
  \cite{matsumura}*{Cor. 29.2} we have $W'(k_R)\cong W(k_R)$; we
  therefore identify $W'(k_R)$ with $W(k_R)$, so $W(k_A)\to
  W(k_R)\subset R$ is $\mf_{W(k_R)}$-\'etale.  Our derivation
  $\partial : W(k_A)\to A$ induces a homomorphism of rings
  \begin{displaymath} \phi(\partial) : W(k_A)\to \frac{A[X]}{(X^2)}
\subset R[X]/(X^2),
\end{displaymath} extending the inclusion $W(k_A)\to R$, and by $\mf_{W(k_R)}$-\'etaleness it lifts uniquely to
a continuous homomorphism
  \begin{displaymath} \phi': W(k_R) \to (\varprojlim
R/\mf_R^n)[X]/(X^2) = \frac{R[X]}{(X^2)},
  \end{displaymath} extending the inclusion $W(k_R)\to R$;
such a homomorphism is of the form $\phi' = \phi(\partial')$
where $\partial'$ is a derivation $W(k_R)\to R$.  This
implies
  \begin{displaymath}
    \partial (W(k_R)) = \partial'(W(k_R))\subset R.
  \end{displaymath} By assumption the local ring extensions
  $A/W(k_A)$,  $R/A$ are tamely
  ramified, implying that  $R/W(k_R)$ is tamely ramified.  The
  extension $R/W(k_R)$ is described by  the relation
  for a uniformising parameter $x\in \mf_R \setminus \mf_R^2$:
  there exists an Eisenstein polynomial $f(X)= X^n +
  a_1X^{n-1} + \cdots + a_n\in W(k_R)[X]$, $p|a_i$, $1\leq
  i\leq n$, $p^2{\not|\ }a_n$ (the divisibility takes place in
  $W(k_R)$), such that $f(x)=0$; hence $p = u x^n$, $u\in
  R\setminus \mf_R$; therefore  by tameness $p{\not|\ } n$.
  Applying $\partial$ gives $f'(x)\partial (x) + f^\partial
  (x)=0$, where $f^\partial (x) = x^{n-1}\partial (a_1) +
  \cdots +\partial (a_n)$.  Since $\partial (a_i)\in R$ and
  $\partial (p)=0$ we have $\partial (a_i) = r_i p = r'_i
  x^n$, for some $r_i, r'_i \in R$, and it follows that one
  can write $f^\partial (x)= rx^n$, for some $r\in R$.  Since
  $p{\not|}\ n$ we have $f'(x) = nx^{n-1} + (n-1)a_1x^{n-2} +
  \cdots + a_{n-1} = u x^{n-1}$, where $u\in R\setminus
  \mf_R$. It follows that $\partial (x)= - f^\partial (x)/
  f'(x) = -rx/u \in \mf_R \subset R$

  (b) $A$ and $R$ are complete of equal characteristic: Then
$A$ has a coefficient field $k_A \subset A$ which injects
into $R$ and this image is part of a coefficient field $k_R$
of $R$, so $k_A \subset k_R \subset R= k_R[[t]]$, where
$t\in \mf_R \setminus \mf_R^2$.  Since $k_R/k_A$ is
separably algebraic there exists to each $\alpha \in k_R $ a
polynomial $f\in k_A[X]$ such that $f(\alpha)=0$,
$f'(\alpha)\neq 0$. Applying $\partial $ gives
$f'(\alpha)\partial (\alpha) + f^\partial (\alpha )=0$,
hence $\partial (\alpha)\in k_R\subset R$. This implies
  \begin{displaymath}
    \partial (k_R)\subset R.
  \end{displaymath} It remains to see that $\partial (t)\in
Rt$. Assume the contrary, $\partial(t)\notin (t)$.  Choose
$x\in A$ such that $R\mf_A = Rx$, and write $x=ut^n$, where
$u$ is a unit in $R$.  Since $\partial (k_R)\subset R$ we
have $\nu (\partial (u)) \geq \nu(\partial(t))$.  Therefore,
since $n$ and $p$ are coprime, $\nu(\partial (x)) < \nu(x)
$.  This contradicts the assumptions $\partial (x)\in \mf_A$
and $\nu(\mf_A)=\nu(x)$.  
%
 
$(3)$: Let $\partial \in T_{A/k}(\log^\pi \mf_A)$. Then for $a\in A$
we have $\partial (a)= ra\in A$ for some $r\in R$.  Any element $z\in
K(A)$ can be written $a/b$ with $a,b\in \mf_A$, implying
  \begin{displaymath}\label{base-relation}
\frac{\partial(z)}z= \frac {\partial (a)}a - \frac{\partial
(b)}b \in R.
\end{displaymath} Since $K(R)=K(A)$ it follows $\partial
(R)\in R$ and $\partial (\mf_{R})\subset \mf_{R}$; moreover
$T_{A/k}^\pi \cap T_{R/k}(\mf_R)\subset T_{A/k}(\mf_A)$. (When
$R/A$ is tamely ramified the inclusion is an equality by $(1)$.)

(4): By $(1)$ $T_{A/k}^\pi \subset T_{R/k}(\mf_R)\cap T_{A/k}(\log^\pi
\mf_A)\subset T_{A/k}(\mf_A) $, so it remains to see that
$T_{A/k}(\log^\pi \mf_A)\subset T_{A/k}^\pi$.  Put $R_1= R\cap K(A)$,
which again is a discrete valuation ring and $k_R/k_{R_1}$ is
algebraic \Lem{\ref{q-finite}}.  We have inclusions of local rings
  \begin{displaymath} A\overset{\pi_1}\hookrightarrow R_1
\overset{i}\hookrightarrow R,
  \end{displaymath} so $\pi = i \circ \pi_1$. That $\pi$ is
tamely ramified implies that $\pi_1$ and $i$ are tamely
ramified.  Let $\partial \in T_{A/k}(\log^\pi \mf_A)$. Then
for $a\in A$ we have $\partial (a)= ra\in A$ for some $r\in
R$, hence $r\in R\cap K(A) = R_1$, hence $\partial (a)\in
R_1a$; therefore
  \begin{displaymath} T_{A/k}(\log^{\pi} \mf_A)\subset
T_{A/k}(\log^{\pi_1} \mf_A)\subset T_{A/k}(\mf_A)\cap
T_{R_1/k}(\mf_{R_1})
\end{displaymath} where the last inclusion follows from
$(3)$.  In particular, $T_{A/k}(\log^{\pi}\mf_A) \subset
T_{R_1/k}(\mf_{R_1})$.  Since $i$ is tamely ramified and
$k_R/k_{R_1}$ is algebraic, by (2) we have
$T_{R_1/k}(\mf_{R_1}) \subset T^i_{R_1/k}$; therefore
$T_{A/k}(\log^\pi \mf_A) \subset T^\pi_{A/k}$.
\end{pfof}

\begin{remark}\label{gabrielov}
If $A\to R$ is a local injective map of
  N{\oe}therian rings, it does not follow that the induced map of
  completed rings $\bar \pi: \bar A \to \bar R$ is injective (see
  \cite{gabrielov,hubl}) but derivations of $A$ extend to derivations
  of $\bar A$ that preserve the kernel of $\bar \pi$, since the
  identity $\pi\circ \partial -
    \partial \circ \pi =0$ extends, by continuity, to the
identity $\bar \pi \circ \partial =
    \partial \circ \bar \pi $. 
\end{remark}
\subsection{Proof of the lifting theorem}

  \begin{thm}\label{mainthm} Let $\pi: X/S \to Y/S$ be a morphism of
    schemes over a scheme $S$.  Assume that $Y$ and $S$ are
    N{\oe}therian and $X$ is a Krull scheme.  Assume that $\pi$ is
    dominant, generically algebraic and tamely ramified
    \Defn{\ref{tame-ramification}}. If $Y$ is of unequal
    characteristic, assume that $Y$ is absolutely tamely ramified. Let
    $D_\pi$ be the discriminant set of $\pi$. Then:
    \begin{enumerate}[label=($L_\theenumi$),ref=($L_\theenumi$)]
    \item Let $J_{D_\pi}$ be {\it any}\/ defining ideal of the closure of
      $D_\pi$. Then
      \begin{displaymath} T_{Y/S}^\pi = T_{Y/S}(\log^\pi J_{D_\pi}).
      \end{displaymath}
      \label{L1}\\
    \item If the residue field extension $k_c/k_{\pi(c)}$ is algebraic
      for each point $c\in C_\pi$ of height $1$, then
      \begin{displaymath} T_{Y/S}^\pi = T_{Y/S}(I_{D_\pi}),\label{L2}
      \end{displaymath}
      where $I_{D_\pi}$ is the (reduced) ideal of the closure of the
      discriminant set $D_\pi$.
    \item \label{L3} Let $I_{C_\pi}$ be the ideal of the critical locus
      of $\pi$.  Then $\pi^{-1}(T_{Y/S}^\pi) \subset T_{X/S}(I_{C_\pi})$, i.e.
      the liftable vector fields are tangent to the critical locus.
    \end{enumerate}
  \end{thm}
  \begin{pf} Since $X/Y$ is generically separably algebraic any
    derivation $\partial $ of $\Oc_Y(U)$, where $U$ is open in $Y$, has a
    unique lift $\bar \partial $ to a derivation of the field $\Oc_{X,\xi}$
    over the generic point $\xi$ of $X$.  Since $X$ is Krull it
    suffices to see that $\bar \partial (\Oc_{X,c}) \subset \Oc_{X,c}$ when $c\in
    X$ is a point of height $1$ that maps to $U$, for then $\bar \partial
    (\Oc_{\pi^{-1}(U)}) \subset \Oc_{\pi^{-1}(U)}$.  If $c\notin C_\pi$, i.e.
    $d\pi_c$ is surjective, this is evidently the case, so assume $c\in
    C_\pi$ is a point of height $1$ in $X$.  Let $c$ be a point of
    $C_\pi$ of height $1$.  Put $R=\Oc_{X,c}$, $A=\Oc_{Y, \pi(c)}$, so
    we may regard $A$ as a subring of the discrete valuation ring $R$,
    with inclusion morphism $\pi : A \to R$.  Put also $k=\Oc_{S,s}$.
    Let $\mf_A$ and $\mf_R$ be the maximal ideal of $A$ and $R$, so
    $\mf_A \subset \mf_R$.  We need to prove that $T^\pi_{A/k} =T_{A/k}
    (\log^\pi \mf_A ) $, and if $k_R/k_A$ is algebraic, $T^\pi_{A/k} =
    T_{A/k}(\mf_A)$.  By \Proposition{\ref{basic-prop}}, $(1)\Rightarrow (3)$,
    the morphism $A\to R$ is ramified (hence tamely ramified).
    Therefore~\ref{L1} follows from \Theorem{~\ref{fund-thm}}, (3),
    and~\ref{L2} follows from \Theorem{~\ref{fund-thm}}, (2). The
    assertion in ~\ref{L3} follows from \Theorem{~\ref{fund-thm}},
    (1).
\end{pf}

\subsection{The normalisation morphism} Let $\pi : \bar X/S \to X/S $
be the normalisation of an integral N{\oe}therian scheme.
\begin{thm}\label{seidenberg}
    \begin{enumerate}[label=(\theenumi)]
  \item $ T_{X/S}(\log^\pi I_D) \subset T^\pi_{X/S}$.
  \item If the normalisation morphism $\pi$ is tamely ramified, then
    $T_{X/S}^\pi = T_{X/S}\subset T_{\bar X/S}(I_C)$ (in particular,
    $\pi$ is weakly submersive).
  \end{enumerate}
\end{thm}

\begin{pf} To alleviate the notation slightly we consider only the
  case $S= \Spec \Zb$; the general relative case is no different. Let
  $x\in \bar X$ , $y=\pi(x)$ and put $A=\Oc_x$, $B=\Oc_y$, where $y$
  lies over a prime number $s$.  We have an inclusion of integral
  domains $A\subset B$ with coinciding fraction fields, $K=K(A)=K(B)$,
  and $B$ is integrally closed in $K$.  We need to prove that a
  derivation $\partial: A\to A$ can be lifted to a derivation $\tilde
  \partial : B \to B$. First, $\partial$ can be extended uniquely to a
  derivation $\tilde \partial: K\to K$.  By the Mori-Nagata theorem
  \cite[Th.  33.10]{nagata:local} the normalisation of a N{\oe}therian
  integral domain is a Krull ring, hence $B$ is the intersection of
  the discrete valuation rings $R$ in $K$ containing $B$. Hence we may
  regard $A$ as a subring of a discrete valuation ring $R$ and it
  suffices to see that $\tilde \partial (R)\subset R$.  Localise $A$
  at the prime $p \in \Spec B$ such that $R= B_p$. Then $A_p \subset
  R$ and $\partial (A_p)\subset A_p$, so one may assume $A= A_p$.  We
  have then an algebraic extension $(A,\mf_A, k_A) \to (R, \mf_R,
  k_R)$ of one-dimensional local rings where $R$ is regular.  Let
  $\nu$ be the normalised discrete valuation of $R$. If $A= R$ there
  is nothing to prove, hence one may assume that $A$ is not regular
  (i.e.\ not normal).

  (1): Let $\partial \in T_{A}(\log^\pi \mf_A)$. The proof
that $\partial (R)\subset R$ is similar to the proof of (3)
in \Theorem{\ref{fund-thm}}. Any element in $K(A)= K(R)$ can
be written $z=a/b $ with $a,b\in \mf_A$. Since $\partial
(a)\in Ra$ for each $a\in \mf_A$, we get
\begin{displaymath} \frac{\partial(z)}z= \frac {\partial (a)}a -
  \frac{\partial (b)}b \in R.
\end{displaymath} In particular, $\partial (R)\subset
R$. This proves $T_{A}(\log^\pi \mf_A) \subset T^\pi_{A}$.

(2): Let $\partial $ be a derivation of $A$, inducing a derivation of $K(A)=
K(R)$, and reasoning as in the proof of \Theorem{\ref{fund-thm}},(2),
it induces a derivation $\bar R \to K(\bar R)$; we have $A\subset R \subset \bar
R$.  Since $\bar R \cap K(R)= R$ it suffices to see that $\partial (\bar R)\subset
\bar R$.  Let $t\in \mf_{\bar R}\setminus \mf_{\bar R}^2$ be a uniformising
parameter.  Let $W(k_R)$ either be a complete $p$-ring such that $\bar
R$ is an Eisenstein extension of $W(k_R) $, or a coefficient field in
$\bar R$, and $W(k_A)$ either be a complete $p$-ring or a coefficient
field in $A$ such that $W(k_A)\subset W(k_R)$. Since $B/A$ is integral, so
the localisation $A_q \subset B_p=R$ is integral, where $p$ is a prime
ideal in $B$ and $q=A\cap p$, by assumption the field extension
$k_R/k_{A_q}$ is separably algebraic. We can assume that $A=A_q$. Now
follow the proof of (2) in \Theorem{\ref{fund-thm}}.  First we get
$\partial(W(k_R)) \subset R$, which is a consequence of the fact that
$W(k_R)/W(k_A)$ is $\mf_{W(k_R)}$-\'etale. Note that this follows
without a priori assuming $\partial(\mf_A)\subset \mf_A$.  Secondly, we get $\partial
(\bar R)\subset \bar R$ if $\partial(t)\in \bar R$.  Finally, by
\Theorem{\ref{fund-thm}}, (1), if $T_A=T^\pi_A$, it follows that $T_A =
T_A(\mf_A) \subset T_R(\mf_R)$.

Assume the contrary, that $\nu(\partial (t))<0$, where $\nu$ is the normalised
valuation of $\bar R$. Let $x_1, \dots , x_n$ be a system of
parameters of $\mf_A$ and put $m_i = \nu(x_i)$, where we may assume
$m_1 \leq m_2 \leq \cdots \leq m_n$.  Since $p {\not \vert \ } m_i$, we get
  \begin{equation}\label{w-one} \nu(\partial (x_i)) =
\nu(x_i)+\nu(\frac{\partial(u_i)}{u_i} + m_i \frac {\partial
(t)}{t}) = m_i + \nu(\partial(t))-1,
\end{equation} where $x_i= u_it^{m_i}$ for some unit
$u_i\in R$. The last equality follows since $\partial
(W(k_R)) \subset R $, implying $\nu(\partial (u))\geq
\nu(\partial(t))$ (since $\nu(\partial(t)) <0$); this holds
both in equal and mixed characteristic. Since
$\partial(x_1)\in A$ and $\nu(\partial(t))<0$, so $0 \leq m_1 + \nu(\partial(t))-1 < m_1$, we get $\nu(\partial(t))= 1-m_1$ and
$m_1 \geq 2$.  To get a contradiction it suffices to see
that $m_1 \vert \ m_i$, $i=1,2, \dots , n$, since then $\Zb
m_1 + \Zb m_2 + \cdots + \Zb m_n$ is a proper ideal of $\Zb$
which would be in opposition to  $K(A)= K(R)$.
Since $m_1 \vert \ m_1$, assume $m_1 \vert \ m_i$ when $i=1,2,
\dots , r-1$ for some integer $r\geq 2$.  We get from
(\ref{w-one}) 
  \begin{displaymath} \nu(\partial (x_r)) = m_r -m_1 =
\sum_{i=1}^{r-1} a_i m_i,
  \end{displaymath}
  where $a_i \in \Nb$; the last equality follows since
  $\nu(\partial(x_r))< m_r$.  Hence $m_1 \vert\  m_r$.
\end{pf}
\begin{cor}\label{1-dim} Let $(A,\mf_A)$ be a local
1-dimensional integral N{\oe}therian ring with a tamely
ramified normalisation.  If there exists a derivation
$\partial $ of $A$ and an element $x\in \mf_A\setminus
\mf_A^2$ such that $\partial (x)\notin \mf_A$, then $A$ is
regular.
\end{cor}
\begin{pf} From the above proof we get $T_A = T_A(\mf_A)$ if $A$ is
  1-dimensional and not normal, implying the assertion.
\end{pf}
\begin{remark} \Corollary{\ref{1-dim}} follows from Zariski's lemma
  when $A$ contains the rational numbers and its completion is reduced
  (see \cite[Corollary to Th. 30.1]{matsumura}).  A more instructive
  argument that $T_{A/k}= T_{A/k}(\mf_A)$ for singular one-dimensional
  $k$-algebras $A$ can be achieved in a special case. Assume that $A$
  is a local ring such that the Jacobian criterion of regularity
  holds.  Then the Jacobian ideal $J$ is preserved by $T_{A/k}$, i.e.
  $T_{A/k}= T_{A/k}(J)$. If $A$ is of characteristic $0$ we have
  $T_{A/k}= T_{A/k}(\sqrt{J})$ (see e.g.  \cite{kallstrom:preserve}),
  hence if $A$ is singular we have a proper $T_{A/k}$-invariant radical
  ideal $\sqrt{J}\subset \mf_A$.  Therefore, if $A$ is a singular
  one-dimensional ring it follows that $T_{A/k}= T_{A/k}(\mf_A)$.  In
  \Corollary{\ref{1-dim}} no assumption is made on the characteristic
  or the validity of the Jacobian criterion. However, for wildly
  ramified normalisations the existence of $\partial\in T_{A/k}$ and
  $t\in \mf_A$ such that $ \partial (t)=1$ does not imply that $A$ is
  regular (\Ex{\ref{wild}}).
\end{remark}
\begin{remark}\label{remark-seidenberg} What Seidenberg actually proved is this.  Let $A$ be an
  integral domain, containing the rational numbers $\Qb$, $\bar A $
  its integral closure and $A'$ the ring of quasi-integral elements
  (so $\bar A = A'$ when $A$ is N{\oe}therian), then $T_{A/\Qb} \subset
  T_{A'/\Qb}$, and if the ring of formal power series $\bar A[[t]]$ is
  integrally closed, then $T_{A/\Qb} \subset T_{\bar A/\Qb}$.  Say that a
  derivation $\partial$ of $A$ is integrable if there exists a Hasse-Schmidt
  derivation of the form $(1, \partial , D_2, \dots , )$.  Note that the
  higher derivations $D_i$ are not uniquely determined by $\partial$; see
  \cite{matsumura:integrable}. The idea in Seidenberg's proof is to
  consider only the integrable derivations, giving rise to
  automorphisms of $A[[t]]$, and the point is that derivations are
  always integrable to higher derivations when $\Qb \subset A$.  On the
  other hand, if $\Qb \nsubseteq A$ the question of integrability is a
  delicate matter.  The only general condition I am aware of that all
  derivations of a $k$-algebra $A$ be integrable is that $A/k$ be
  smooth \cite{matsumura:integrable}, but this is of no interest for
  the lifting problem, at least when $k$ is a perfect field, since
  then $A$ is normal. Letting $T_{X/S}^{int}$ be the sheaf of
  integrable derivations, Seidenberg's theorem states that
  $T_{X/S}^{int}\subset T^\pi_{X/S}$ when $\pi$ is the normalisation morphism
  of a N{\oe}therian scheme
  \cite{seidenberg:derivationsintegralclosure} (also when $\pi$ is
  wildly ramified).  Assume that $X/k$ is a variety. Matsumura
  \cite{matsumura:integrable} proves that if $k$ is a perfect field,
  then $T_{X/k}^{int}$ is of the same rank as $T_{X/k}$, and gives an
  example that over an imperfect field $k$ one may have
  $T^{int}_{X/k}=0$ while $T_{X/k}\neq 0$. On the other hand,
  $T_{X/S}(\log^\pi \mf_A)\subset T^\pi_{X/k}$ \Th{\ref{seidenberg}}, and
  this module is always non-zero when $T_{X/k}\neq 0$.
\end{remark}
We extend an example by Seidenberg
\cite{seidenberg:derivationsintegralclosure} giving many examples of
rings with wild normalisations and where there exist derivations that
do not lift.
\begin{example}\label{wild} Let $R$ be a discrete valuation
ring of characteristic $p>0$, $K$ its fraction field and
$\partial $ a derivation of $R$ such that
$\partial(\mf_R)\not \subset \mf_R$; we also regard
$\partial$ as a derivation of $K$ and let $K^\partial\subset
K$ be the subfield of elements $r$ such that $\partial(r)=0$
(so $R^p \subset K^p \subset K^\partial$).  Suppose $S$ is a
subset of $R$ such that $\partial (S)\subset S$ and
$\partial^n(S)=0$ for some positive integer $n$.  The set
$\mf_R^{\partial }= K^\partial \cap \mf_R$ is non-trivial if
and only of $\Char R =p >0$, and then $\mf_R^\partial
\subset \mf_R^p$. Select an element $c\in \mf_R^\partial $
and an increasing sequence of subsets $B_0 \subset B_1
\subset \cdots \subset B_n \subset \mf_R^\partial$ such that
$B_i \subset cB_{i+1}$ and let $A$ be the smallest local
subring of $R$ that contains the set
  \begin{displaymath} S_B=SB_0 + \partial (S)B_1 + \cdots +
\partial^{n-1}(S)B_n.
  \end{displaymath} Since clearly $\frac 1c
\partial(S_B)\subset S_B$ we get $\frac 1c
  \partial (A)\subset A$, but $\frac 1c \partial
(R)\not\subset R$.  If $S$ generates $K$ and $S\cap B_0\neq
\emptyset$, then $K(A)=K$ and the inclusion $\pi: A\to R$ is
the normalisation of $A$. The derivation $\delta = \frac 1c
\partial$ of $A$ does not belong to $ T_{A}(\log^\pi
\mf_A)$.  Let $k$ be a field of characteristic $p$,
$R=k[t]_0$, $S=\{t^p, t^{p+1}\}\subset k[t]_0$ and
$B_i=t^{(i+1)p}$.  Then $A=k[t^p, t^{p+1}]_0 \subset R=
k[t]_{0} \subset k(t)$, which is an extension with
ramification index $p$. This is close to Seidenberg's
example.  Then $\delta =t^{-p}\partial_t$ is a derivation of
$A$ such that $\delta (R)\not \subset R$.  Clearly,
$\delta(t^{p+1})\notin Rt^{p+1}$, so $\delta \notin
T_{A}(\log^\pi \mf_A)$.  Note also that $\partial_t \in
T^\pi_{A}$, but $\partial_t (\mf_R)\not\subset \mf_R$,
showing that $(1)$ in \Theorem{\ref{fund-thm}} need not hold
when $\pi$ is wildly ramified.  
\end{example} All derivations may lift even if the  
normalisation is wildly ramified.
\begin{example}\label{wild-norm} Put $A= k[t^3, t^7,t^8 ]_0 \subset R=
  k[t]_0$ and assume $\Char k = 3$. Here $T_{A/k} = A \partial_t \subset
  T_{R/k}$.
\end{example} 

\section{Birational morphisms}\label{birational}
\subsection{Differentially ramified morphisms} Assuming the morphism
$\pi$ is of the type in \Theorem{\ref{mainthm}} we have $T_{X/Y}=0$.
Therefore, if $J$ is an ideal of $\Oc_X$ it makes sense to define the
sub-sheaf of $T_{X/S}$  of liftable derivations that preserve $J$ by
\begin{displaymath} T^\pi_{Y/S}(J) = T^\pi_{Y/S}\cap T_{X/S}(J).
\end{displaymath} 
For an ideal $I$ of $\Oc_Y$, its inverse image ideal is $\tilde I =
\Oc_XI$ (with obvious interpretation).

\begin{definition} A dominant morphism of schemes $\pi: X\to Y$, where
  $X$ is Krull and $Y$ locally N{\oe}therian, with discriminant locus
  $D$, is {\it uniformly ramified}\/ at a generic point $\xi $ of the
  critical locus $C$ if the stalk $I_{D, \pi(\xi)}$ of the ideal of $D$
  at $\pi(\xi)$ has a basis $I_{D, \pi(\xi)} = (x_1, \dots , x_l)$
  satisfying $\nu_\xi (x_i)= \nu_\xi (x_j)$, for each $ i, j$. The
  morphism $\pi$ is uniformly ramified if it is uniformly ramified at
  each generic point of the critical locus $C$.
\end{definition} Thus the condition is on the base $Y$ but local at
each generic point of $C$. Do not confuse this with another notion of
uniform ramification in \cite{griffith:ram}.

\begin{definition} A generically finite morphism $\pi$ is
{\it differentially ramified}\/ if
  \begin{displaymath} T_{Y/S}^\pi =T_{Y/S}(I_{D}),
  \end{displaymath}
  where $I_D$ is the (reduced) ideal of of the discriminant $D$.
\end{definition} Thus for differentially ramified morphisms
the liftable derivations can be characterised in terms of
the discriminant $D$ only, and do not depend on the
particular differentially ramified morphism with such a
discriminant. For example, residually algebraic morphisms as
in \Theorem{\ref{mainthm}} are differentially ramified.
\begin{thm}\label{blowup1} Let $\pi: X\to Y$ be a morphism as in
  \Theorem{\ref{mainthm}}.
    \begin{enumerate}[label=(\theenumi)]
    \item If $\pi$ is uniformly ramified, then it is differentially
      ramified.
  \item Assume that for each generic point $\xi$ of $C$,
there exists a basis $I_{D,\pi(\xi)}= (x_1, \dots , x_r)$
and derivations $\partial_1, \dots,
    \partial_r \in T_{Y/S, \pi(\xi)}$ such that $\det
\partial_j(x_i)$ is invertible. Then $\pi$ is uniformly
ramified if and only if it is differentially ramified.
  \end{enumerate}
\end{thm} For example, by (1), if $D_\pi$ is a hypersurface,
i.e. $I_{D_\pi}$ is locally principal, then morphisms of the
considered type are differentially ramified.
\begin{pf} (1): If $\pi$ is uniformly ramified a derivation
$\partial \in T_{Y/S}(I_{D})$ satisfies, for each generic
point $\xi$ of $C$, $\nu_\xi (\partial(x_i))\geq \nu_\xi
(x_i)$ for some basis $(x_1, \dots , x_l)$ of $I_{D,
\pi(\xi)}$. The assertion therefore follows from
\Theorem{\ref{mainthm}}.

(2): By (1) we have only to prove that $\pi$ is uniformly ramified when
$T_{Y/S}^\pi= T_{Y/S}(I_{D})$. Let $c_{ij}$ be the inverse matrix of
$(\partial_i (x_j))$. Setting $\partial'_i = \sum c_{ij}\partial_j$ we have $\partial'_i (x_j)=
\delta_{ij}$, so that we can assume $\partial_i(x_j)= \delta_{ij}$.  By assumption $
T_{Y/S}(I_{D})\subset T_{Y/S}(\log^\pi I_{D})$ \Th{\ref{mainthm}}; hence
$\nu_\xi (x_i) = \nu_\xi(x_i \partial_j (x_j)) \geq \nu_{\xi} (x_j)$. By symmetry
the opposite inequality also holds.
\end{pf}
\subsection{Blow-ups} A projective birational morphism of integral
N{\oe}therian schemes $\pi : \tilde X\to X$ can be constructed as a
blow-up of a coherent fractional ideal on $X$, and if $X$ and $Y$ are
projective over a field one may even construct $\pi$ as a blow-up of an
ideal \cite{EGA3}*{Prop. 2.3.5, Cor. 2.3.7}.  Let $X$ be an integral
scheme with sheaf of rational functions $R_X$, let $I$ be a coherent
$\Oc_X$-submodule of $R_X$, and let $\pi: Bl_I(X)\to X$ be the blow-up
of $X$ along $I$. We will in this section try to describe the liftable
derivations $T^\pi_{X/S}$, without using valuations along the critical
locus of $\pi$; more precisely, we look for a (fractional) ideal $J$
such that $T^\pi_{X/S} = T_{X/S}(J)$. Clearly, if such an ideal exists,
then $ D_\pi = \{x\in X \ \vert \ J_x \neq \Oc_{X,x}\}$, the discriminant
locus of $\pi$, so $J$ and $I$ would  have the same radical.

Let $p :\overline {Bl_I(X)} \to Bl_I(X)$ be the normalisation morphism
of the blow-up of $I$. The integral closure of $I$ is the sheaf $\bar
I$ whose sections over an open set $U$ is
\begin{displaymath}
  \bar I (U) = \bigcap_{x\in U, \hto
    (x)\leq 1} \Oc_{\overline{Bl_I(X)},x} I_x.
\end{displaymath}
\begin{lemma}(see
  \cite{vasconcelos:integral_closure})\label{int-clos-lemma}
  \begin{eqnarray*}
  \bar I(U)&=& \{f \in R_X(U) \ \vert \ \text{ there exists an equation } f^n + b_1f^{n-1}+ \cdots + b_n=0,\\ && n>0, b_i\in I^i(U) \}
  \\ &=& (\pi \circ p)_* ( I \Oc_{\overline {Bl_I(X)}}).
\end{eqnarray*}
\end{lemma}
\Lemma{\ref{int-clos-lemma}} implies that $\bar I$ is coherent 
since $\pi \circ p$ is proper. Clearly, $I\subset \bar I$, and if $X$ is
normal and $I$ an ideal of $\Oc_X$ then $\bar I\subset \Oc_X$. See
\cite{henzer-johnston.lantz.shah:ratliff-rush} for relations between
$\overline {Bl_I(X)}$ and $Bl_{\overline I^n}(X)$ for high $n$.

Define $[I^{n+1}: I^n]= \{a \in R_X \ \vert \ aI^n \subset I^{n+1}\}$; this
is a coherent fractional ideal.  Then for sufficiently high $n$ we
have $[I^{n+2}: I^{n+1}]= [I^{n+1}: I^n]$.  Assuming that $I$ be
coherent, the Ratliff-Rush fractional ideal $\hat I$ associated to $I$
is the common fractional ideal $ [I^{n+1}: I^n]$ for high $n$, and we
have $\hat I \subset \bar I$.  One has $I^n = (\hat I)^n$ for high $n$
[loc. cit.]. These assertions are well-known when $I$ is a coherent
ideal \cite{ratliff-rush} and they immediately generalise to
fractional ideals since a coherent fractional ideal is locally
contained in an ideal of the form $ J/f$, where $J$ is a coherent
ideal of $\Oc_X$, and we have $[I^{n+1}:I^n]= \frac 1f [(fI)^{n+1}:
(fI)^n]$.

\begin{thm}\label{blow-up-thm} Let $X/S$ be an integral N{\oe}therian scheme with sheaf of
  rational functions $R_X$. Let $I$ be a coherent $\Oc_X$-submodule of
  $R_X$ and $\pi : Bl_I(X)\to X$ be the blow-up of $X$ over $I$.
  \begin{enumerate}[label=($\theenumi$),ref=($\theenumi$)]
  \item  $T_{X/S}(\hat I) \subset T^\pi_{X/S} \subset T_{X/S}(\bar I)$ (the right
    inclusion holds when the normalisation morphism $p$ is tamely
    ramified).\label{no1} \\
  \item  $T_{X/S}(I)\subset T^\pi_{X/S}\cap T_{X/S}(\bar I)\cap T_{X/S}(\hat I)$. \label{no2}
  \end{enumerate}
  Assume that $p$ is tamely ramified.  If $\hat I = \bar I$ (e.g. $I=
  \bar I$), then $T^\pi_{X/S}= T_{X/S}(\bar I) $. If $I$ is radical,
  then $\pi$ is differentially ramified.
\end{thm}

\begin{pf} We first prove $T_{X/S}(I)\subset T^\pi_{X/S}$.  A derivation
  that preserves $I$ extends to a derivation of the Rees algebra
  $\oplus_{n\geq 0} I^n$, and therefore induces a derivation of each affine
  chart $B_{(f)}= \cup_{j\geq 0} I^j/f^j$, $f \in I$, of $Bl_I(X)=\Proj
  \oplus_{n\geq 0} I^n$.  Since in such a chart $\tilde I = BI (= fB)$,
  clearly $\partial (\tilde I )\subset \tilde I$.

  $ T^\pi_{X/S} \subset T_{X/S}(\bar I)$: We have $T^\pi_{X/S} \subset T^{\pi \circ
    p}_{X/S}$ since $p$ is tamely ramified \Th{\ref{seidenberg}}.  As
  $\overline{Bl_I(X)}$ is Krull, by \Theorem{\ref{mainthm}}, ($L_3$),
  a liftable derivation $\partial $ lifts to a derivation $\tilde \partial \in
  T_{\overline {Bl_I(X)}/S}$ that preserves the ideal of the critical
  divisor. Since $\tilde I =I\Oc_{\overline {Bl_I(X)}}$ is a locally
  principal ideal defining the critical locus of $\pi \circ p$ it follows
  that $\tilde \partial ( \tilde I)\subset \tilde I$. Therefore $\partial$ preserves
  $(\pi \circ p)_*(\tilde I) = \bar I$.

  $ T_{X/S}(\hat I) \subset T^\pi_{X/S} $: Let $\hat \pi : Bl_{\hat I}(X)\to
  X$ be the blow-up of $\hat I$. We then have
\begin{displaymath}
  T_{X/S}(\hat I)\subset T^{\hat  \pi}_{X/S} =  T^{\pi}_{X/S}
\end{displaymath}
where the latter equality follows since $\hat I ^n = I^n$ when $n\gg
1$, so $Bl_{\hat I}(X)= Bl_I(X)$. This completes the proof of the
assertions \ref{no1} and \ref{no2}.

Assume $\hat I = \bar I$, which holds in particular when $I$ i
radical, then \ref{no1}  implies $T^\pi_{X/S}= T_{X/S}(\bar I)$, which in
particular implies that $\pi$ is differentially ramified when $I$ is
radical.
\end{pf}

\begin{cor}\label{cor:blow} Let $\pi: \tilde X/S \to X/S$ be a tamely ramified blow-up
  of a reduced subscheme $V/S$ of an integral N{\oe}therian scheme
  $X/S$, such that $\tilde X$ is Krull. If $X/S$ is smooth at generic
  points of $V$, then $\pi$ is uniformly ramified.
\end{cor}
Suppose $s$ is a point in $S$ with perfect residue field and that the
fibre $V_s$ is a smooth subscheme of $X_s$, then $\tilde X_s$ is
regular and hence Krull (normal).
\begin{pf} For each generic point of the critical locus its image $x$
  is a generic point of $ V \subset X$. Let $k$ be the residue field of the
  image $s$ of $x$ in $S$.  The canonical exact sequence $0 \to
  \mf_{X,x}/\mf_{X,x}^2\to k_{X,x}\otimes \Omega_{X/S,x}\to \Omega_{k_{X,x}/k} \to 0$
  is (split) exact since $k_{X,x}/k$ is formally smooth. Assume that
  $\{x_1, \dots , x_r\}\subset \mf_{X,x}$ induces a basis $\{\bar x_1,
  \dots , \bar x_r\}$ of $\mf_{X,x}/\mf_{X,x}^2$, so $\mf_{X,x}=(x_1,
  \dots , x_r)$, and select elements $\omega_i\in \Omega_{X/S,x}$ such that
  $\{x_1, \dots , x_r, \omega_1, \dots , \omega_s\}$ induces a basis $\{1\otimes
  d_{X/S}(x_1), \dots , 1\otimes d_{X/S}(x_r), 1\otimes \omega_1 , \dots , 1\otimes
  \omega_s\}$ of $k_{X,x}\otimes \Omega_{X/S,x}$. Since $\Omega_{X/S, x}$ is a free
  $\Oc_{X,x}$-module of finite type it follows that
  $\{d_{X/S}(x_1),\dots , d_{X/S}(x_r), \omega_1, \dots , \omega_s\}$ are free
  generators of $\Omega_{X/S,x}$; hence there exist derivations $\partial_i \in
  T_{X/S,x}$ satisfying $\partial_i(x_j)=\delta_{ij}$.  Since $V$ is reduced and
  $\pi$ is tamely ramified it follows that $\pi$ is differentially
  ramified \Th{\ref{blow-up-thm}}, so the assertion can be concluded
  from \Theorem{\ref{blowup1}}, (2).
\end{pf}
For an ideal $I$ of a polynomial ring $R$ over a field it is in
general a laborious task to compute its integral closure $\bar I$
\cite{vasconcelos:integral_closure}, although for monomial ideals it
is simply the monomial ideal defined by the convex hull of the
exponent vectors of $I$. The Ratliff-Rush ideal associated with $I$
seems even more difficult to compute (see
\cite{rossi-swanson:ratliff-rush,elias:ratliff-rush}); for instance, I
am unaware of any simple description of $\hat I$ for monomial $I$.
Notice that to compute the liftable derivations $T^\pi_{X/S}$ for a
blow-up using \Theorem{\ref{mainthm}} we need to know that $Bl_I(X)$
is normal. This however is difficult to read off from $I$ (a
sufficient condition is that $I^n$ be integrally closed for high $n$).
\begin{example} 
\label{rr-closure-ex}
Let $R=\Qb[[t^4, t^5, t^6,t^7]]$ (singular ring) and $I=(t^4, t^5)$.
Then $I^2=(t^8,t^9,t^{10})$, $[I^2:I]=I$, $[I^3:I^2] = (t^4, t^5,
t^6,t^7)= \mf = \hat I$. We have $T_{R/\Qb}(I) = R t\partial_t$ and
$T_{R/\Qb}(\hat I) = Rt\partial t + R t^2\partial_t$. Hence $T_{R/\Qb}(\hat I)\neq
T_{R/\Qb}(I)$.
    \end{example}
      \begin{example}\label{examples-der2}
        Let $R=\Qb[x,y]_{(x,y)}$ (regular ring), $I= (x^{10}, x^8y,
        xy^4, y^5) $. Then $\hat I = (x^{10},y^5, xy^4$
        $,x^7y^2,x^6y^3,x^8y)$
        (\cite{rossi-swanson:ratliff-rush,elias:ratliff-rush}) and
        $\bar I = (x^{10},x^8y,x^6y^2, x^4y^3, xy^4,y^5)$.  Generators
        of the module of derivations of a monomial $\mf$-primary ideal
        $J= (x^{a_i}y^{b_i})_{i=0}^r $, $a_i > a_{i+1}$ are easily
        computed.  Let $h= \max\ (a_{i}- a_{i+1}) $ and $w= \max \
        (b_{i+1}-b_i)$.  Then $T_R(J) = R x\partial_x + Ry\partial_y + Ry^{h}\partial_x
        + R x^w\partial_y$. In particular, $T_{R/\Qb}(I) = (x,y^3)\partial_x +
        (x^7,y)\partial_y$, $T_{R/\Qb}(\hat I)= (x,y)\partial_x + (x^5,y)\partial_y$,
        and $T_{R/\Qb}(\bar I) = (x,y)\partial_x +(x^3,y)\partial_y$, so
        $T_{R/\Qb}(I)\neq T_{R/\Qb}(\hat I)\neq T_{R/\Qb}(\bar I)$.
        Consider the morphism
        \begin{displaymath}
          \pi : Bl_I(\Spec R) \to \Spec R.
        \end{displaymath}
        Since $\hat I \neq \bar I$ it follows that $Bl_I(\Spec R)$
        cannot be normal, so we cannot use \Theorem{\ref{mainthm}} to
        compute $T_{R/\Qb}^\pi$; we know only that $(x,y)\partial_x +
        (x^5,y)\partial_y \subset T^\pi_{R/\Qb}\subset (x,y)\partial_x + (x^3,y)\partial_y $. It
        would require some effort to determine whether the derivations
        $x^3\partial_y, x^4\partial_y$ belong to $T^\pi_{R/\Qb}$.
      \end{example}

\begin{remark}\label{examples-der}
  We have $T_{X/S}(I)\subset T_{X/S}(I^n)\subset T_{X/S}([I^n: I^{n-1}])$, but
  in general $T_{X/S}(I^n)\not\subset T_{X/S}(I)$, when $I \neq \hat I$, as
  seen from \Example{\ref{rr-closure-ex}},  for $n=3$.
  \end{remark} 
  We have $T_{X/S}(I)\subset T_{X/S}(\hat I)\subset T_{X/S}(\bar I)$ where the
  first inclusion follows from \Remark{\ref{examples-der}}, while the
  latter, valid when $p$ is tamely ramified, follows from
  \Theorem{\ref{blow-up-thm}}; that these inclusions in general are
  strict is shown in \Example{\ref{examples-der2}}. One can refine
  these inclusions.  In \cite{shah:coefficient} Shah defines for
  $\mf$-primary ideals in a local ring the notion of coefficient
  ideals, using coefficients of the Hilbert polynomial of $I$; this
  notion was extended in
  \cite{henzer-johnston.lantz.shah:ratliff-rush}*{Def 3.21} to general
  ideals $I$ in a N{\oe}therian domain, defining coefficient ideals as
  contractions of pull-backs of $I$ by certain birational morphisms.
  Let $I$ be a coherent ideal of $\Oc_X$, where $X$ is an integral and
  N{\oe}therian scheme of dimension $d$, and let $I_{\{k\}}$, $k=0,
  \dots , d$ be the coefficient ideals of $I$ as defined in [loc.
  cit.]; put also $I_{\{d+1\}}=I$.  One has $I \subset I_{\{d\}} = \hat I
  \subset I_{\{d-1\}}\subset \cdots \subset I_{\{1\}} \subset I_{\{0\}}=\bar I$ [loc.  cit]
  and one can prove, as in the beginning of the proof of
  \Theorem{\ref{blow-up-thm}}, noting that the lifted derivations
  evidently preserve the intersections of local rings described in
  [loc. cit, 3.14], that $T_{X/S}(I) \subset T_{X/S}(I_{\{d\}} )=
  T_{X/S}(\hat I) \subset T_{X/S}(I_{\{d-1\}})\subset \cdots\subset T_{X/S}(I_{\{k\}})
  \subset\cdots \subset T_{X/S}(I_{\{1\}}) \subset T_{X/S}(\bar I)$, at least if the
  birational map that is used to define $I_{\{k\}}$ is tamely
  ramified.  Let $J_{\{k\}}$ be the smallest coherent ideal that
  contains $I_{\{k\}}$ such that $T_{X/S}(I_{\{k-1\}}) \subset
  T_{X/S}(J_{\{k\}})$.  We then have
  \begin{displaymath}
    I_{\{k\}}\subset J_{\{k\}} \subset I_{\{k-1\}},
  \end{displaymath}
  and one may wonder when these inclusions are strict? In
  \Examples{\ref{examples-der}} and {\ref{examples-der2}} we have
  $J_{\{d+1\}}= I_{\{d\}} (= \hat I)$.  If one knows a priori that
  $J_{\{k\}} = I_{\{k-1\}}$, one can compute $I_{\{k-1\}}$ from a
  knowledge of $I_{\{k\}}$ and $T_{X/S}(I_{\{k-1\}})$, since
  $J_{\{k\}}=\Dc_{X/S}(T_{X/S}(I_{\{k-1\}}))I_{\{k\}} $ and the ring
  $\Dc_{X/S}(T_{X/S}(I_{\{k-1\}}))$ of differential operators
  generated by $T_{X/S}(I_{\{k-1\}})$. However, one cannot use
  \Theorem{\ref{mainthm}} to describe $T_{X/S}(I_{\{k-1\}})$, since
  the source of the birational morphisms used to define $I_{\{k-1\}}$
  is not Krull.

  We end this paper with two concrete illustrations of the use of
  \Theorems{\ref{fund-thm}}{\ref{blow-up-thm}}.
\begin{example} Let $A$ be the polynomial ring $k[y_1, \dots
, y_r]$ localised at some prime ideal, and assume that $A
\subset R$ is tamely ramified and $K(R)/K(A)$ is finite,
where $R$ is a discrete valuation ring with valuation $\nu$.
\Theorem{\ref{fund-thm}} implies $\partial =
a_1\partial_{y_1} + \cdots + a_r\partial_{y_r}\in
T_{A/k}^\pi$, $a_i\in A$, if and only if $\nu(a_i)\geq
\nu(y_i)$, $i=1, \dots , r$. To be very explicit, let $A=
k[y_1,y_2]_{(y_1,y_2)}$, $R= k[x_1,x_2]_{(x_1)}$, so $\mf_A=
(y_1,y_2)$ and $\mf_R=(x_1)$, and define a birational map
$\pi: A \to R$ by $y_1 = x_1$, $y_2= x_1^nx_2$ where $n$ is
an integer $\geq 1$. Then $\partial_{y_1}=
  \partial_{x_1} -n x_2/x_1\partial_{x_2}$ and
$\partial_{y_2} = x_1^{-n}\partial_{x_2} $, acting on $K(R)=
k(x_1, x_2)$.  The condition $\partial \in T_{A/k}(\log^\pi
\mf_A)$ is that $\partial \in T_{A/k}$ and $\nu(\partial
(y_1)) \geq \nu(y_1)= 1$, $\nu (\partial (y_2))\geq
\nu(y_2)$, i.e. $a_1\in (y_1,y_2)$ and $a_2\in (y_1^n,y_2)$;
therefore $T_{A/k}(\log^\pi \mf_A) =\mf_A
  \partial_{y_1} + (y_1^n, y_2) \partial_{y_2}$.  If $\Char
k =0$ or $n$ is coprime to $\Char k$, by
\Theorem{\ref{fund-thm}}
  \begin{displaymath} T^\pi_{A/k} = T_{A/k}(\log^\pi \mf_A)
= \mf_A \partial_{y_1} + (y_1^n, y_2)\partial_{y_2}.
  \end{displaymath} We have then
  \begin{displaymath} RT_{A/k}^\pi = R\nabla_1 +
R\partial_{x_2} = T_{R/k}(I_{C_\pi}),
\end{displaymath} where the critical ideal $I_{C_\pi}$ is
$(x_1)$. If $n=1$, so $\pi$  comes by a chart of  the blow-up  of $\mf_A$,
then $T^\pi_{A/k} =T_{A/k}(\log^\pi \mf_A)=
\mf_A\partial_{y_1} + \mf_A\partial_{y_2}= T_{A/k}(\mf_A)$, so $\pi$ is differentially ramified,  in agreement with \Theorem{\ref{blow-up-thm}}, noting that $\mf_A$ is radical.  If $n\geq
2$, corresponding to   a successive blowing up of points, then the
derivations $y_1^k\partial_{y_2}\in T_{A/k}(\mf_A)$, $k=
1,\dots , n-1$ are not liftable. In this case $A \to R$ comes by a chart of the
blow-up of the ideal $I=(y_1^n,y_2)$; this ideal is Ratliff-Rush
and integrally closed $I = \hat I = \bar I$. It is clear that $T_{A/k}^\pi = \mf_A
\partial_{y_1} + (y_1^n, y_2)\partial_{y_2} = T_{A/k}(I)$, as predicted by \Theorem{\ref{blow-up-thm}}. If $ \Char
k =p>0$ and $p$ divides $n$, so $\pi$ has wild ramification, then
$T^\pi_{A/k}= A \partial_{y_1} + (y_1^n,y_2)
\partial_{y_2}$. 
\end{example}
\begin{example} Put $B =\Qb[x,s,t]/(s^4-xt)$ and $A= \Qb[x,y]$ and
  define the homomorphism $\pi : A\to B$, $\pi(x)=x, \pi(y)= x^2s$. The
  critical divisor is given by the principal ideal $(s)\subset B$ and $x=
  s^4/t$ in $B_{(s)}$, so the corresponding valuation $\nu: A \to \Zb$
  is determined by $\nu(x)= 4, \nu(y)= 9$. A derivation $a\partial_x + b\partial_y\in
  T_{A/\Qb}$ lifts to a derivation of $B$ if $\nu(a)\geq \nu(x)= 4$ and
  $\nu(b)\geq \nu(y)= 9$ \Th{\ref{fund-thm}}, noting that $B$ is a normal
  ring, implying $T^\pi_{A/\Qb} = (x,y)\partial_x + (x^2,y)\partial_y$.
\end{example}


\begin{bibsection}
  \begin{biblist} 
\bib{abhyankar:res-surf}{book}{
  author={Abhyankar, S. S.},
  title={Resolution of singularities of embedded algebraic surfaces},
  series={Springer Monographs in Mathematics},
  edition={2},
  publisher={Springer-Verlag},
  place={Berlin},
  date={1998},
  pages={xii+312},
  isbn={3-540-63719-2},
  review={MR1617523 (99c:14021)},
}

\bib{arnold:wave}{article}{
  author={Arnol{\cprime }d, V. I.},
  title={Wave front evolution and equivariant Morse lemma},
  journal={Comm. Pure Appl. Math.},
  volume={29},
  date={1976},
  number={6},
  pages={557\ndash 582},
  review={MR 55 \#9148},
}

\bib{auslander-buchsbaum}{article}{
  author={Auslander, M.},
  author={Buchsbaum, D. A.},
  title={On ramification theory in noetherian rings},
  journal={Amer. J. Math.},
  volume={81},
  date={1959},
  pages={749\ndash 765},
  review={MR 21 \#5659},
}

\bib{bierstone-milman:resolution}{article}{
  author={Bierstone, Edward},
  author={Milman, Pierre D.},
  title={Canonical desingularization in characteristic zero by blowing up the maximum strata of a local invariant},
  journal={Invent. Math.},
  volume={128},
  date={1997},
  number={2},
  pages={207\ndash 302},
  issn={0020-9910},
  review={MR1440306 (98e:14010)},
}

\bib{borel:Dmod}{book}{
  author={Borel, A},
  author={Grivel, P.-P},
  author={Kaup, B},
  author={Haefliger, A},
  author={Malgrange, B},
  author={Ehlers, F},
  title={Algebraic $D$-modules},
  series={Perspectives in Mathematics},
  publisher={Academic Press Inc.},
  address={Boston, MA},
  date={1987},
  volume={2},
}

\bib{bourbaki:commutative}{book}{
  author={Bourbaki, Nicolas},
  title={Commutative algebra},
  publisher={Hermann},
  date={1972},
}

\bib{elias:ratliff-rush}{article}{
  author={Elias, Juan},
  title={On the computation of the Ratliff-Rush closure},
  journal={J. Symbolic Comput.},
  volume={37},
  date={2004},
  number={6},
  pages={717\ndash 725},
  issn={0747-7171},
  review={MR2095368 (2005j:13022)},
}

\bib{encinas-villamayor:good}{article}{
  author={Encinas, S.},
  author={Villamayor, O.},
  title={Good points and constructive resolution of singularities},
  journal={Acta Math.},
  volume={181},
  date={1998},
  number={1},
  pages={109\ndash 158},
  issn={0001-5962},
  review={MR1654779 (99i:14020)},
}

\bib{fossum:krull}{book}{
  author={Fossum, Robert M.},
  title={The divisor class group of a Krull domain},
  note={Ergebnisse der Mathematik und ihrer Grenzgebiete, Band 74},
  publisher={Springer-Verlag},
  place={New York},
  date={1973},
  pages={viii+148},
  review={MR0382254 (52 \#3139)},
}

\bib{gabrielov}{article}{
  author={Gabri{\`e}lov, A. M.},
  title={Formal relations among analytic functions},
  language={Russian},
  journal={Izv. Akad. Nauk SSSR Ser. Mat.},
  volume={37},
  date={1973},
  pages={1056\ndash 1090},
  issn={0373-2436},
  review={MR0346184 (49 \#10910)},
}

\bib{griffith:ram}{article}{
  author={Griffith, Phillip},
  title={Some results in local rings on ramification in low codimension},
  journal={J. Algebra},
  volume={137},
  date={1991},
  number={2},
  pages={473\ndash 490},
  issn={0021-8693},
  review={MR1094253 (92c:13017)},
}

\bib{EGA1}{article}{
  author={Grothendieck, A.},
  title={\'El\'ements de g\'eom\'etrie alg\'ebrique. I. Le langage des sch\'emas},
  journal={Inst. Hautes \'Etudes Sci. Publ. Math.},
  number={4},
  date={1960},
  pages={228},
  issn={0073-8301},
  review={MR0217083 (36 \#177a)},
}

\bib{EGA4:1}{article}{
  author={Grothendieck, A.},
  title={\'El\'ements de g\'eom\'etrie alg\'ebrique. IV. \'Etude locale des sch\'emas et des morphismes de sch\'emas. I},
  language={French},
  journal={Inst. Hautes \'Etudes Sci. Publ. Math.},
  number={20},
  date={1964},
  pages={259},
  issn={0073-8301},
  review={MR0173675 (30 \#3885)},
}

\bib{EGA4:4}{article}{
  author={Grothendieck, A},
  title={\'El\'ements de g\'eom\'etrie alg\'ebrique. IV. \'Etude locale des sch\'emas et des morphismes de sch\'emas IV},
  language={French},
  journal={Inst. Hautes \'Etudes Sci. Publ. Math.},
  number={32},
  date={1967},
  pages={361},
}

\bib{EGA3}{article}{
  author={Grothendieck, Alexander},
  author={Dieudonne, Jean},
  title={\'{E}tude cohomolique des faisceaux coh{\'e}rent},
  date={1961},
  journal={Publ. IHES},
  number={4},
}

\bib{henzer-johnston.lantz.shah:ratliff-rush}{article}{
  author={Heinzer, William},
  author={Johnston, Bernard},
  author={Lantz, David},
  author={Shah, Kishor},
  title={Coefficient ideals in and blowups of a commutative Noetherian domain},
  journal={J. Algebra},
  volume={162},
  date={1993},
  number={2},
  pages={355\ndash 391},
  issn={0021-8693},
  review={MR1254782 (95b:13020)},
}

\bib{hubl}{article}{
  author={H{\"u}bl, Reinhold},
  title={Completions of local morphisms and valuations},
  journal={Math. Z.},
  volume={236},
  date={2001},
  number={1},
  pages={201\ndash 214},
  issn={0025-5874},
  review={MR1812456 (2002d:13003)},
}

\bib{kallstrom:purity}{unpublished}{
  author={K{\"a}llstr{\"o}m, Rolf},
  title={Purity of branch, critical, and discriminant locus},
  date={2006},
  note={Submitted, avaialble at arXiv},
}

\bib{kallstrom:preserve}{article}{
  author={K{\"a}llstr{\"o}m, Rolf},
  title={Preservation of defect sub-schemes by the action of the tangent sheaf},
  journal={J. Pure and Applied Algebra},
  volume={156},
  date={2005},
  number={2},
  pages={286\ndash 319},
  issn={0001-8708},
  review={MR 2001m:58078},
}

\bib{lazarsfeld:vol2}{book}{
  author={Lazarsfeld, Robert},
  title={Positivity in algebraic geometry. II},
  series={Ergebnisse der Mathematik und ihrer Grenzgebiete. 3. Folge. A Series of Modern Surveys in Mathematics [Results in Mathematics and Related Areas. 3rd Series. A Series of Modern Surveys in Mathematics]},
  volume={49},
  note={Positivity for vector bundles, and multiplier ideals},
  publisher={Springer-Verlag},
  place={Berlin},
  date={2004},
  pages={xviii+385},
  isbn={3-540-22534-X},
  review={\MR {2095472 (2005k:14001b)}},
}

\bib{lipman:equisingular}{article}{
  author={Lipman, Joseph},
  title={Equisingularity and simultaneous resolution of singularities},
  conference={ title={Resolution of singularities}, address={Obergurgl}, date={1997}, },
  book={ series={Progr. Math.}, volume={181}, publisher={Birkh\"auser}, place={Basel}, },
  date={2000},
  pages={485--505},
  review={\MR {1748631 (2001e:32044)}},
}

\bib{matsumura:integrable}{article}{
  author={Matsumura, Hideyuki},
  title={Integrable derivations},
  journal={Nagoya Math. J.},
  volume={87},
  date={1982},
  pages={227\ndash 245},
}

\bib{matsumura}{book}{
  author={Matsumura, Hideyuki},
  title={Commutative ring theory},
  publisher={Cambridge University Press},
  date={1986},
}

\bib{nagata:local}{book}{
  author={Nagata, Masayoshi},
  title={Local rings},
  note={Corrected reprint},
  publisher={Robert E. Krieger Publishing Co., Huntington, N.Y.},
  date={1975},
  pages={xiii+234},
  isbn={0-88275-228-6},
  review={MR0460307 (57 \#301)},
}

\bib{neron:minimaux}{article}{
  author={N{\'e}ron, Andr{\'e}},
  title={Mod\`eles minimaux des vari\'et\'es ab\'eliennes sur les corps locaux et globaux},
  language={French},
  journal={Inst. Hautes \'Etudes Sci. Publ.Math. No.},
  volume={21},
  date={1964},
  pages={128},
  issn={0073-8301},
  review={\MR {0179172 (31 \#3423)}},
}

\bib{popescu:neron}{article}{
  author={Popescu, Dorin},
  title={General N\'eron desingularization and approximation},
  journal={Nagoya Math. J.},
  volume={104},
  date={1986},
  pages={85--115},
  issn={0027-7630},
  review={\MR {868439 (88a:14007)}},
}

\bib{ratliff-rush}{article}{
  author={Ratliff, L. J., Jr.},
  author={Rush, David E.},
  title={Two notes on reductions of ideals},
  journal={Indiana Univ. Math. J.},
  volume={27},
  date={1978},
  number={6},
  pages={929\ndash 934},
  issn={0022-2518},
  review={MR0506202 (58 \#22034)},
}

\bib{rossi-swanson:ratliff-rush}{article}{
  author={Rossi, Maria Evelina},
  author={Swanson, Irena},
  title={Notes on the behavior of the Ratliff-Rush filtration},
  booktitle={Commutative algebra (Grenoble/Lyon, 2001)},
  series={Contemp. Math.},
  volume={331},
  pages={313\ndash 328},
  publisher={Amer. Math. Soc.},
  place={Providence, RI},
  date={2003},
  review={MR2013172 (2005b:13006)},
}

\bib{scheja:fortzetsungderivationen}{article}{
  author={Scheja, G\"unther},
  author={Storch, Uwe},
  title={Fortsetzung von Derivationen},
  language={German},
  journal={J. Algebra},
  volume={54},
  date={1978},
  number={2},
  pages={353\ndash 365},
}

\bib{seidenberg:derivationsintegralclosure}{article}{
  author={Seidenberg, A},
  title={Derivations and integral closure},
  journal={Pacific J. Math.},
  volume={16},
  date={1966},
  pages={167\ndash 173},
}

\bib{shah:coefficient}{article}{
  author={Shah, Kishor},
  title={Coefficient ideals},
  journal={Trans. Amer. Math. Soc.},
  volume={327},
  date={1991},
  number={1},
  pages={373--384},
  issn={0002-9947},
  review={\MR {1013338 (91m:13008)}},
}

\bib{spivakovsky:popescu}{article}{
  author={Spivakovsky, Mark},
  title={A new proof of D. Popescu's theorem on smoothing of ring homomorphisms},
  journal={J. Amer. Math. Soc.},
  volume={12},
  date={1999},
  number={2},
  pages={381\ndash 444},
  issn={0894-0347},
  review={MR1647069 (99j:13008)},
}

\bib{vasconcelos:integral_closure}{book}{
  author={Vasconcelos, Wolmer},
  title={Integral closure},
  series={Springer Monographs in Mathematics},
  note={Rees algebras, multiplicities, algorithms},
  publisher={Springer-Verlag},
  place={Berlin},
  date={2005},
  pages={xii+519},
  isbn={978-3-540-25540-6},
  isbn={3-540-25540-0},
  review={\MR {2153889 (2006m:13007)}},
}

\bib{villamayor:equisingular}{article}{
  author={Villamayor U., Orlando},
  title={On equiresolution and a question of Zariski},
  journal={Acta Math.},
  volume={185},
  date={2000},
  number={1},
  pages={123--159},
  issn={0001-5962},
  review={\MR {1794188 (2002a:14003)}},
}

\bib{villamayor:res}{article}{
  author={Villamayor, Orlando},
  title={Constructiveness of Hironaka's resolution},
  journal={Ann. Sci. \'Ecole Norm. Sup. (4)},
  volume={22},
  date={1989},
  number={1},
  pages={1\ndash 32},
  issn={0012-9593},
  review={MR985852 (90b:14014)},
}

\bib{zariski:equisinglarI}{article}{
  author={Zariski, Oscar},
  title={Studies in equisingularity. I. Equivalent singularities of plane algebroid curves},
  date={1965},
  journal={Amer. J. Math.},
  volume={87},
  pages={507\ndash 536},
  review={ MR 31 \#2243},
}


  \end{biblist}
\end{bibsection}
\end{document}